\documentclass{amsart}
\usepackage[all]{xy}
\usepackage{color}
\usepackage{amsthm}% allows theorem* , e.g.
\numberwithin{equation}{subsection}% makes equat numb contain the section

\swapnumbers %pour que le numéro apparaisse devant le théorème

\newtheorem{theorem}[equation]{Theorem}% note that equations and thms get the same numbering: makes more sense.
\newtheorem*{theorem*}{Theorem}
\newtheorem{lemma}[equation]{Lemma}

\newtheorem{proposition}[equation]{Proposition}
\newtheorem*{proposition*}{Proposition}
\newtheorem{corollary}[equation]{Corollary}
\newtheorem*{corollary*}{Corollary}

\theoremstyle{remark}
\newtheorem{definition}[equation]{Definition}
\newtheorem*{definition*}{Definition}

\newtheorem{example*}{Example}

\newtheorem{notation}[equation]{Notation}

\theoremstyle{remark}
\newtheorem{remark}[equation]{Remark}
\newcommand{\op}{^{\mathsf{op}}}
\newcommand{\cA}{{\mathcal A}}
\newcommand{\cAo}{{\mathcal A}\op}
\newcommand{\cCo}{{\mathcal C}\op}
\newcommand{\cB}{{\mathcal B}}
\newcommand{\cBo}{{\mathcal B}\op}
\newcommand{\cBh}{\widehat{\cB}}
\newcommand{\cAh}{\widehat{\cA}}
\newcommand{\cC}{{\mathcal C}}
\newcommand{\cD}{{\mathcal D}}

\newcommand{\cI}{{\mathcal I}}

\newcommand{\cM}{{\mathcal M}}

\newcommand{\cT}{{\mathcal T}}

\newcommand{\cU}{{\mathcal U}}

\newcommand{\bbL}{\mathbb{L}}

\newcommand{\bbU}{\mathbb{U}}
\newcommand{\bbV}{\mathbb{V}}

\newcommand{\bbS}{\mathbb{S}}
\newcommand{\bbW}{\mathbb{W}}
\newcommand{\bbZ}{\mathbb{Z}}

\newcommand{\Ab}{\mathsf{Ab}}
\newcommand{\Spe}{\mathsf{Spt}}% stabilization
\newcommand{\rep}{\mathsf{rep}}
\newcommand{\orep}{\overline{\mathsf{rep}}}
\newcommand{\A}{\mathsf{A}}
\newcommand{\Amod}{\cA\text{-}\mbox{Mod}}
\newcommand{\Aomod}{\cA\op\text{-}\mbox{Mod}}

\newcommand{\TAomod}{T(\cA)\op\text{-}\mbox{Mod}}

\newcommand{\Bmod}{\cB\text{-}\mbox{Mod}}
\newcommand{\Bomod}{\cB\op\text{-}\mbox{Mod}}

\newcommand{\TBomod}{\cT(\cB)\op\text{-}\mbox{Mod}}

\newcommand{\Sp}{\mathsf{Sp}^{\Sigma}}%Symmetric spectra
\newcommand{\Ape}{\cA^{\sharp}_{pe}}
\newcommand{\hm}{\underline{h}}
\newcommand{\per}{\mathsf{per}}
\newcommand{\Add}{\mathsf{Add}}
\newcommand{\Bpe}{\cB^{\sharp}_{pe}}

\newcommand{\Dsma}{\wedge^{\bbL}}
\newcommand{\Spcat}{\mathsf{Sp}^{\Sigma}\text{-}\mbox{Cat}}%Spectral categories
\newcommand{\Ho}{\mathsf{Ho}}

\newcommand{\too}{\longrightarrow}

\newcommand{\ie}{\textsl{i.e.}\ }

\begin{document}

\title[Matrix invariants of spectral categories]{Matrix invariants of spectral categories}
\author{Gon{\c c}alo~Tabuada}

\address{Gon\c calo Tabuada, Departamento de Matem{\'a}tica e CMA, FCT-UNL, Quinta da Torre, 2829-516 Caparica,~Portugal }
\email{tabuada@fct.unl.pt}
\date{\today}


\begin{abstract}
In this paper, we further the study of spectral categories, initiated in~\cite{Spectral}. Our main contribution is the construction of the {\em Universal matrix invariant of spectral categories}, \ie a functor $\cU$ with values in an additive category, which inverts the Morita equivalences, satisfies matrix invariance, and is universal with respect to these two properties. For example, the algebraic $K$-theory and the topological Hochschild and cyclic homologies are matrix invariants, and so they factor uniquely through $\cU$. As an application, we obtain for free non-trivial trace maps from the Grothendieck group to the topological Hochschild homology ones.
\end{abstract}

\maketitle

\tableofcontents
%-----------------------------------------------------------------------
\section{Introduction}
%-----------------------------------------------------------------------
\subsection*{Spectral categories} {\em Spectral categories} are categories enriched over the symmetric monoidal category $\Sp$ of symmetric spectra \cite{HSS}.  As linear categories can be understood as {\em rings with several objects}, spectral categories can be understood as {\em symmetric ring spectra with several objects}. The precise statement is that a symmetric ring spectra is a spectral category with a single object. Due to this ``flexibility'', spectral categories pervade several mathematical areas: Blumberg-Mandell's work~\cite{Mandell} on topological Hochschild homology and its variants; Schwede-Shipley's work~\cite{SS} on the classification of stable model categories;  Kontsevich's non-commutative algebraic geometry program~\cite{Kontsevich} \cite{finMotiv}; Dugger's work \cite{Dugger} on spectral enrichments of model categories;  $\ldots$ 
In this article, we further the study of spectral categories initiated in~\cite{Spectral}, as follows: \subsection*{Matrix invariants}
All the classical (functorial) invariants such as algebraic $K$-theory $K(-)$, topological Hochschild homology $THH(-)$, topological cyclic homology~$TC(-)$, $\ldots$ extend naturally from symmetric ring spectra to spectral categories. See Chapters \ref{sec:AlgK}-\ref{sec:THH} for details. In a ``motivic spirit'' (see Kontsevich's talk~\cite{Kontsevich}), we would like to study all these classical invariants simultaneously. By analysing the commum features of all them, we were lead to the following notion of {\em matrix invariant}. Some definitions are in order. As in the case of ring spectra, given a spectral category $\cA$, we can consider its derived category $\cD(\cAo)$ of right $\cA$-modules (\S\ref{sub:review}). A spectral functor $\cA \to \cB$ is called a {\em Morita equivalence} (\ref{def:Morita}) if its restriction of scalars functor $\cD(\cBo) \stackrel{\sim}{\to} \cD(\cAo)$ is an equivalence. An {\em upper triangular matrix}  $\underline{M}$ (\ref{UTM}) is given by
$$
\begin{array}{rcl}
\underline{M} & := & \begin{pmatrix} \cA & X \\ \ast &
  \cC \end{pmatrix}\,,
\end{array}
$$
where $\cA$ and $\cC$ are spectral categories and $X$
is a $\cA\mbox{-}\cC$-bimodule. The totalization $|\underline{M}|$ of $\underline{M}$ is the spectral category whose set of objects is the disjoint union of the sets of objects of $\cA$
and $\cC$ and whose morphisms are given by: $\cA(x,y)$ if $x, y \in \cA$; $\cC(x,y)$ if $x, y \in \cC$ ; $X(x,y)$ if $x \in \cA,\, y \in \cC$ and $\ast$ if $x \in \cC,\, y \in \cA$. The composition is induced by the composition on $\cA$, $\cC$ and the $\cA\text{-}\cC$-bimodule structure on $X$. Notice that we have two natural inclusion spectral functors
$$
\begin{array}{lccr}
i_1: \cA \too |\underline{M}| && i_2: \cC \too |\underline{M}| \,.
\end{array}
$$
\begin{definition*}
Let $F: \Spcat \to \A$ be a functor, from the category of (small) spectral categories, to an additive category $\A$. We say that $F$ is a {\em matrix invariant of spectral categories} if it verifies the following two conditions:
\begin{itemize}
\item[M)] it sends the Morita equivalences to isomorphisms in $\A$ and
\item[MA)] for every upper triangular matrix $\underline{M}$, the inclusion spectral functors induce an isomorphism in $\A$
$$ [F(i_1) \,\, F(i_2)]: F(\cA) \oplus F(\cC) \stackrel{\sim}{\too} F(|\underline{M}|)\,.$$
\end{itemize}
\end{definition*}
In Propositions~\ref{prop:AddKth} and \ref{prop:THHAdd}, we prove that $K(-)$, $THH(-)$, $TC(-)$ are all examples of matrix invariants of spectral categories. Now a natural question~arises:

\vspace{0.1cm}

{\it Question (1): Does there exists a ``Universal'' matrix invariant $\cU:\Spcat \to \Add$ of spectral categories ? In the case of existence, and since universal objects tend to be rather formal, can we describe $\cU$ and $\Add$ explicitly ? }

\vspace{0.1cm}

In order to solve {\it Question (1)}, we need first to describe the localization of $\Spcat$ with respect to the class of Morita equivalences. In \cite{Spectral}, we have constructed a Quillen model structure on $\Spcat$, whose weak equivalences are the stable quasi-equivalences (a particular class of Morita equivalences, see~\ref{def:quasi}). Let us denote by $\Ho(\Spcat)$ the homotopy category hence obtained. 

Given spectral categories $\cA$ and $\cB$, we describe in Theorem~\ref{thm:Hom}, the Hom-set $[\cA, \cB]$ in $\Ho(\Spcat)$ in terms of isomorphism classes of a certain derived category of $\cA\mbox{-}\cB$-bimodules $\rep(\cA,\cB)$. This answers affirmatively to a question raised by To{\"e}n in \cite{Toen}. An important ingredient for this description is a careful adaptation of Lurie's path object construction~\cite[A.3.4.8]{Lurie} to the case of spectral categories, see Theorem~\ref{thm:path}. 

In Chapter~\ref{sec:chapter3} we introduce the notion of {\em spectral triangulated category}. Roughly, it consists of a ``spectral enhancement'' of the classical notion of (idempotent complete) triangulated category~\cite{Neeman}. See~\ref{def:triang} for details. Using Theorem~\ref{thm:Hom} and some general arguments developped by To{\"e}n~\cite{Toen}, we prove in Theorem~\ref{thm:envelope}, that the inclusion $\Ho(\Spcat)^{tr}\subset \Ho(\Spcat)$ of the spectral triangulated categories admits a left adjoint functor $(-)^{\sharp}_{pe}$, which we refer to as the {\em triangulated envelope}.
Our first main Theorem is the following:
\begin{theorem*}{(\ref{thm:localization})}
The composition
$$ \Spcat \too \Ho(\Spcat) \stackrel{(-)^{\sharp}_{pe}}{\too} \Ho(\Spcat)^{tr} $$
is the localization functor associated to the class of Morita equivalences.
\end{theorem*}
The next step towards the {\em ``Universal''} matrix invariant, is the additivization of the category $\Ho(\Spcat)^{tr}$. We proceed as follows: given $\cA, \cB \in \Ho(\Spcat)^{tr}$, the category $\rep(\cA, \cB)$ carries a natural triangulated structure and so we can consider its Grothendieck group $K_0\rep(\cA,\cB)$. Let $\Add$ be the category whose objects are those of $\Ho(\Spcat)^{tr}$ and whose Hom-sets are defined as
$$ \Add(\cA, \cB):= K_0\rep(\cA,\cB)\,.$$
The composition is the induced one. We have a functor $\Ho(\Spcat)^{tr} \to~\Add$, which sends each isomorphism class of $\rep(\cA, \cB)$ to the corresponding class in the Grothendieck group $K_0\rep(\cA,\cB)$. Our second main Theorem is the solution to {\em Question (1)}:
\begin{theorem*}{(\ref{thm:main})}
The composed functor
$$\cU: \Spcat \too \Ho(\Spcat) \stackrel{(-)^{\sharp}_{pe}}{\too} \Ho(\Spcat)^{tr} \too \Add$$
is the Universal matrix invariant of spectral categories, \ie for every additive category $\A$, the functor $\cU$ induces a bijection between the additive functors from $\Add$ to $\A$ and the matrix invariants of spectral categories with values in $\A$. 
\end{theorem*}
At this point another natural question arises:

\vspace{0.1cm}

{\it Question (2): According to the preceding Theorem, the additive category $\Add$ contains all the ``information'' about the classical matrix invariants. Hence, how can we decode this information ?}

\subsection*{Grothendieck group and trace maps}
Let $\underline{\bbS}$ be the spectral category with a single object and endomorphisms the sphere symmetric ring spectrum $\bbS$.
Our (partial) solution to {\em Question (2)} is the following:
\begin{proposition*}{(\ref{prop:co-repres})}
For every spectral category $\cB$, we have a natural isomorphism of abelian groups
$$ \Add(\cU(\underline{\bbS}), \cU(\cB)) \stackrel{\sim}{\too} K_0(\cB)\,.$$
\end{proposition*}
This co-representability result has the following important application: 
\begin{corollary*}{(\ref{cor:trace})}
Let $j$ be a non-negative integer and
$$ 
\begin{array}{lcr}
K_0(-):\Spcat \too \Ab && THH_j(-) : \Spcat \too \Ab\,,
\end{array}
$$
the Grothendieck and the $j$-th topological Hochschild homology group functors. Then, each generator $g$ of the $j$-th stable homotopy group of the sphere, furnishes us for free a non-trivial trace map
$$ tr_{j,g}:K_0(-) \Rightarrow THH_j(-)\,.$$
\end{corollary*}

Let us exemplify the preceding Corollary by recalling from \cite[I-Example 2.1]{Schwede} the stable homotopy groups of the sphere up to dimension $8$:

\vspace{0.1cm}

\begin{center}
  \begin{tabular}{ c | c | c | c | c | c | c | c | c | c | c }
    j & 0 & 1 & 2 & 3 & 4 & 5 & 6 & 7 & 8 & $\cdots$ \\ \hline
    $\pi_j^s\bbS$ & $\bbZ$ & $\bbZ/2$ & $\bbZ/2$ & $\bbZ/24$ & 0 & 0 & $\bbZ/2$ & $\bbZ/240$ & $(\bbZ/2)^2$ & $\cdots$ \\ \hline
    $\mbox{generator}$ $g$ & $\iota$ & $\eta$ & $\eta^2$ & $\nu$ &  &  & $\nu^2$ & $\sigma$ & $\eta\sigma, \epsilon$ & $\cdots$ \\
  \end{tabular}
\end{center}

\vspace{0.1cm}

Blumberg and Mandell proved in~\cite[Thm.\,1.3]{Mandell} that given a quasi-compact and semi-separated scheme $X$, its topological Hochschild homology as defined by Geisser and Hesselholt in \cite{GH}, can be recovered from the topological Hochschild homology of the spectral category ``naturally'' associated to the dg category of perfect
complexes on $X$. See \cite{SpectralAlg} for the precise relationship between (the homotopy theories of) spectral and dg categories. Therefore, an application of Corollary~\ref{cor:trace} to this algebraic geometric setting furnishes for free non-trivial trace maps.
\subsection*{Related works}
The ``motivic'' idea of constructing universal invariants is not new and appears in several different subjects: Voevodsky's work~\cite{Voevodsky} on algebraic geometry; Meyer-Nest's work~\cite{Nest} on $C^{\ast}$-algebras; Corti{\~n}as-Thom's work~\cite{Cortinas} on bivariant $K$-theory; Garkusha's work~\cite{Garkusha} on associative rings, the author's work~\cite{IMRN}~\cite{IMRNC}~\cite{Duke} on dg categories, $\dots$ The work presented here is morally the ``topological version'' of the differential graded case \cite{IMRN}, where the notion of matrix invariance corresponds to the notion of additivity. However, we would like to emphasize that since the category of symmetric spectra is not additive (but only up to homotopy), the main constructions and key arguments from \cite{IMRN} are not available. For instance, the analogue of Theorem~\ref{thm:localization} is proved in \cite{IMRN} using model structures constructed using the additivity of the category of (co)chain complexes, wherein here we work at the homotopic level (see the detailed explanations above). To the best of the author's knowledge, the results in this paper offer the first ``motivic'' construction in the context of brave new algebra~\cite{EKMM}\cite{HSS}\cite{Lydakis}\cite{May}.

\vspace{0.2cm}

\noindent\textbf{Acknowledgments\,:} It is a great pleasure to thank Stefan Schwede and Bertrand To{\"e}n for motivating conversations and Gustavo Granja for comments on an older version of this article. 

%-----------------------------------------------------------------------
\section{Preliminaries}\label{sec1}
%-----------------------------------------------------------------------

\subsection{Notations}
Let $\cM$ be a Quillen model category~\cite{Quillen}. We will denote by $[-,-]$ the Hom-sets in its homotopy category $\Ho(\cM)$. Let $\Sp$ be the category of symmetric spectra (of pointed simplicial sets)~\cite{HSS} \cite{Schwede}, endowed with its  projective stable model structure~\cite[III-Thm.\,2.2]{Schwede}. Recall that its set of generating cofibrations consists of
$$I_{proj}=\{ F_n\partial \Delta[m]^+ \too F_n\Delta[m]^+\}_{n,m\geq 0} \,\,,$$
where $F_n(-)$ is the $n$-th free symmetric spectrum construction~\cite[I-2.12]{Schwede}. We denote by $\ast$ the initial and terminal object in $\Sp$, by $-\wedge-$ the smash product bi-functor and by $\mathbb{S}$ its unit, \ie the sphere symmetric spectrum. See \cite[I-\S3]{Schwede} for details. Throughout the article the adjunctions will be displayed vertically, with the left, resp. right, adjoint on the left-hand side, resp. right-hand side.

\vspace{0.1cm}

The results in this article are most conveniently stated for small (spectral) categories. However, in Chapters \ref{sec:1}-\ref{sec:chapter3}  we will need ``large'' (spectral) categories. Therefore, we will use the language of Grothendieck universes~\cite{Grothendieck} to make our statements and proofs rigorous. In what follows $\bbU \in \bbV  \in \bbW \ldots$ will denote Grothendieck universes. When the universes are irrelevant, we will omit them.
 
%\subsection{Triangulated categories}{(Expandir esta seccao-Neeman???)}
%\begin{definition}{(\cite[4.2.7]{Neeman})}\label{def:compact}
%Let $\cT$ be a triangulated category with infinite sums. A full triangulated subcategory of $\cT$ is called {\em localizing} if it is closed under coproducts. A set $\cP$ of objects of $\cT$ is called a set of {\em generators} if the onl localizing subcategory which contains the objects of $\cP$ is $\cT$ itself. An object $X \in \cT$ is {\em compact} if for any family of objects $\{A_i\}_{i \in I}$, the canonical map
%$$ \underset{i \in I}{\bigoplus} \,\cC(X, A_i) \too \cC(X, \underset{i \in I}{\bigoplus} A_i)$$
%is an isomorphism.
%\end{definition}
%\begin{proposition}{(\cite[?]{Neeman})}\label{prop:Neeman}
%Let $F: \cT_1 \to \cT_2$ be a triangulated functor which commutes with infinite sums. Suppose that $\cT_1$ and $\cT_2$ are compactly generated and there is a set $\cG\subset \cT_1$ of compact generators which is stable under (co)suspensions. Consider the following conditions:
%\begin{itemize}
%\item[(a)] For all $X, Y \in \cG$, the triangulated functor $F$ induces a bijection
%$$ \cT_1(X,Y) \too \cT_2(FX, FY)\,.$$
%\item[(b)] The set of objects $\{FX |\, X \in \cG \}$ is a set of compact generators in $\cT_2$.
%\end{itemize}
%If $F$ satisfies condition (a), then $F$ is fully faithful. If moreover $F$ satisfies condition (b), it is an equivalence of triangulated categories.
%\end{proposition}

\subsection{Review on spectral categories}\label{sub:review}

References on spectral categories are \cite[\S 2]{Mandell}, \cite[Appendix\,A]{SS} and \cite[\S 2]{Spectral}. Recall that a {\em $\bbU$-small spectral category $\cA$} consists of the following data: 
\begin{itemize}
\item[-] a $\bbU$-small set of objects $\mbox{obj}(\cA)$ (usually denoted by $\cA$ itself); 
\item[-] for each pair of objects $(x,y)$ of $\cA$, a $\bbU$-small symmetric spectrum $\cA(x,y)$; 
\item[-] for each triple of objects $(x,y,z)$ of $\cA$, a composition morphism in $\Sp_{\bbU}$
$$\cA(y,z)\wedge \cA(x,y) \too \cA(x,z)\,,$$
satisfying the usual associativity condition; 
\item[-] for any object $x$ of $\cA$, a morphism $\bbS \to \cA(x,x)$ in $\Sp_{\bbU}$, satisfying the usual unit condition with respect to the above composition.
\end{itemize}
We denote by $\Spcat_{\bbU}$ the category of $\bbU$-small spectral categories. 
Let $\cA$ be a (fixed) $\bbU$-small spectral category. A {\em $\bbU$-small $\cA$-module} is a morphism $\cA \to \Sp_{\bbU}$ in $\Spcat_{\bbV}$. We denote by $\Amod_{\bbU}$ the category of $\bbU$-small $\cA$-modules. By \cite[Thm.\,A.1.1]{SS}, $\Amod_{\bbU}$ is a $\bbV$-small (cofibrantly generated) $\Sp_{\bbU}$-model category (\ref{def:SpeMod}). We denote by $\cD(\cA)$ the {\em derived category} of $\cA$, \ie the homotopy category $\Ho(\Amod_{\bbU})$. Notice that we have two (fully faithful) {\em Yoneda morphisms} in $\Spcat_{\bbV}$
$$
\begin{array}{lccr}
\begin{array}{rcl}
\hm^{-}: \cA\op & \to & \Amod_{\bbU}\\
z & \mapsto & \cA(z,-)
\end{array}
& &
\begin{array}{rcl}
\hm_{-}: \cA & \to &\Aomod_{\bbU}\\
z & \mapsto & \cA(-,z)
\end{array}\,,
\end{array}
$$
where $\cA\op$ is the {\em opposite spectral category of $\cA$}, \ie $\cA\op$ has the same objects as $\cA$ and $\cA\op(x,y)=\cA(y,x)$. By \cite[\S A.1]{SS}, a morphism $F:\cA \to \cB$ in $\Spcat_{\bbU}$ gives rise to a {\em restriction/extension of scalars} spectral Quillen adjunction (on the left)
$$\xymatrix{
\Bomod_{\bbU} \ar@<1ex>[d]^{F^{\ast}} &&& \cD(\cBo) \ar@<1ex>[d]^{F^{\ast}} \\
\Aomod_{\bbU} \ar@<1ex>[u]^{F_!} &&& \cD(\cAo) \ar@<1ex>[u]^{\bbL F_!}\,,
}
$$
which can be naturally derived (on the right).
\subsection{Quillen model structure}
Given a spectral category $\cA$, we can form a genuine category $[\cA]$ by keeping the same set of objects and defining the set of morphisms between $x$ and $y$ in $[\cA]$ to be the set of morphisms $[\bbS, \cA(x,y)]$ in $\Ho(\Sp)$. We obtain in this way a functor
$$ [-]: \Spcat \too \mbox{Cat}\,,$$
with values in the category of small categories.
\begin{definition}\label{def:quasi}
A spectral functor $F:\cA \to \cB$ is a {\em stable quasi-equivalence} if:
\begin{itemize}
\item[S1)] for all objects $x,y \in \cA$, the morphism in $\Sp$
$$ F(x,y):\cA(x,y) \too \cB(Fx,Fy) $$
is a stable equivalence~\cite[II-4.1]{Schwede} and
\item[S2)] the induced functor
$$ [F]:[\cA] \too [\cB]$$
is an equivalence.
\end{itemize}
\end{definition}
Notice that if $F$ satisfies condition S1), then condition S2) is equivalent to:
\begin{itemize}
\item[S2')] the induced functor
$$ [F]:[\cA] \too [\cB]$$
is essentially surjective.
\end{itemize}
\begin{theorem}{(\cite[Thm.\,5.10]{Spectral})}\label{thm:Model}
The category $\Spcat$ admits a right proper Quillen model structure whose weak equivalences are the stable quasi-equivalences.
\end{theorem}
We denote by $\Ho(\Spcat)$ the corresponding homotopy category. We obtain then an induced functor
$$[-]: \Ho(\Spcat) \too \Ho(\mbox{Cat})\,,$$
with values in the category of small categories and isomorphism classes of functors between them. Moreover, by construction, the natural functor
$$ \Ho(\Spcat_{\bbU}) \too \Ho(\Spcat_{\bbV})$$
is fully faithful.
\begin{proposition}{(\cite[Prop.\,5.13]{Spectral})}\label{prop:fibrant}
A spectral category $\cA$ is fibrant, with respect to the model structure of Theorem~\ref{thm:Model}, if and only if for all objects $x,y \in \cA$ the symmetric spectrum $\cA(x,y)$ is an $\Omega$-spectrum.
\end{proposition}
\begin{proposition}{(\cite[Prop.\,4.18]{Spectral})}\label{prop:cof}
Let $\cA$ be a cofibrant spectral category, with respect to the model structure of Theorem~\ref{thm:Model}. Then for all objects $x,y \in \cA$ the symmetric spectra $\cA(x,y)$ is cofibrant.
\end{proposition}
\begin{remark}\label{rk:cof}
By construction of the generating cofibrations \cite[4.4]{Spectral}, there exists a cofibrant replacement functor $Q(-)$ on $\Spcat$, such that for any spectral category $\cA$, the natural spectral functor $Q(\cA) \to \cA$ induces the identity map on the sets of objects.
\end{remark}
%--------------------------------------------------------------------
\section{Generalized (bi)modules}\label{sec:1}
%--------------------------------------------------------------------
\begin{definition}{(\cite[3.5.1]{SS})}\label{def:SpeMod}
A {\em $\bbV$-small $\Sp_{\bbU}$-model category} is a $\bbV$-small model category $\cM$ which is tensored, cotensored and enriched (denoted by $\underline{\cM}$) over $\Sp_{\bbU}$, such that the following compatibility axiom holds:
\begin{itemize}
\item[(SP)] for every cofibration $A \rightarrow B$ and every fibration $X \rightarrow Y$ in $\cM$, the induced map
$$ \underline{\cM}(B,X) \longrightarrow \underline{\cM}(A,X) \times_{\underline{\cM}(A,Y)} \underline{\cM}(B,Y)$$
is a stable projective fibration in $\Sp_{\bbU}$. If in addition one of the maps is a stable equivalence, then the resulting map is also a stable equivalence. We use the notation $K \wedge X$ and $X^K$ to denote the tensors and cotensors for $X \in \cM$ and $K \in \Sp_{\bbU}$.
\end{itemize}
Given a $\bbV$-small $\Sp_{\bbU}$-model category $\cM$, its {\em internal spectral category} $\mbox{Int}(\cM) \subset \underline{\cM}$ consists on the full $\bbV$-small spectral subcategory of fibrant and cofibrant objects.
\end{definition}
\begin{remark}\label{rk:chunck}
Since every object in $\mbox{Int}(\cM)$ is fibrant and cofibrant, the compatibility axiom (SP) implies that for all objects $X, Y \in \mbox{Int}(\cM)$, the symmetric spectrum $\mbox{Int}(\cM)(X,Y)$ is an $\Omega$-spectrum. Therefore, by Proposition~\ref{prop:fibrant}, $\mbox{Int}(\cM)$ is a $\bbV$-small fibrant spectral category. Notice that we have an equivalence of $\bbV$-small categories
$[\mbox{Int}(\cM)] \simeq \Ho(\cM)$. Moreover, $\mbox{Int}(\cM)$ satisfies the following condition: given two morphisms $X \to Y$ and $X \to Y'$ in $\mbox{Int}(\cM)$, there exists a (functorial) factorization of $X \to Y\times Y'$ in $\cM$:
$$ X \stackrel{p}{\too} \overline{X} \stackrel{q}{\too} Y \times Y'\,,$$ 
where $p$ is a trivial cofibration, $q$ is a fibration and $\overline{X}\in \mbox{Int}(\cM)$.  
\end{remark}
\begin{notation}\label{not:hat}
Given a $\bbU$-small fibrant spectral category $\cB$, we will denote by $\cBh$ the $\bbV$-small spectral category $\mbox{Int}(\Bomod_{\bbU})$.
\end{notation}

\subsection{Generalized modules}
The notion of module over a spectral category admits the following generalization: let $\cM$ be a $\bbV$-small cofibrantly generated $\Sp_{\bbU}$-model category. Then for every $\bbU$-small spectral category $\cA$, we can consider the $\bbV$-small category $\cM^{\cA}$ of morphisms in $\Spcat_{\bbV}$ from $\cA$ to $\underline{\cM}$. The category $\cM^{\cA}$ is endowed with a (cofibrantly generated) Quillen model structure for which the weak equivalences and fibrations are defined objectwise. Moreover the $\Sp_{\bbU}$-enrichment of $\cM$ endow $\cM^{\cA}$ with a natural structure of $\Sp_{\bbU}$-model category. Notice that when $\cM=\Sp_{\bbU}$, we recover the notion of $\bbU$-small $\cA$-module. Finally, a morphism $F:\cA \to \cB$ in $\Spcat_{\bbU}$ gives rise to a {\em restriction/extension of scalars} spectral Quillen adjunction (on the left)
$$\xymatrix{
\cM^{\cB} \ar@<1ex>[d]^{F^{\ast}} &&& \Ho(\cM^{\cB}) \ar@<1ex>[d]^{F^{\ast}} \\
\cM^{\cA} \ar@<1ex>[u]^{F_!} &&& \Ho(\cM^{\cA}) \ar@<1ex>[u]^{\bbL F_!}\,,
}
$$
which can be naturally derived (on the right). 
\begin{proposition}{(\cite[Prop.\,3.2]{Toen})}\label{prop:Toen1}
Let $F:\cA \to \cB$ be a stable quasi-equivalence between $\bbU$-small spectral categories and $\cM$ a $\bbV$-small cofibrantly generated $\Sp_{\bbU}$-model category. 
Assume that the domains of the generating cofibrations of $\cM$ are cofibrant and that for every cofibrant object $X \in \cM$, and every stable equivalence $Z \to Z'$ in $\Sp_{\bbU}$, the induced map (\ref{def:SpeMod})
$$ Z \wedge X \too Z'\wedge X$$
is a weak equivalence in $\cM$. Then the Quillen adjunction $(F_!,F^{\ast})$ is a Quillen equivalence.
\end{proposition}
\begin{proof}
The proof is analogous to the one of \cite[Prop.\,3.2]{Toen}. Replace the tensor product $-\otimes-$ by the smash product $-\wedge-$ and the notion of $C(k)$-model category (see \cite[\S 3]{Toen}) by the notion of $\Sp_{\bbU}$-model category (\ref{def:SpeMod}).
\end{proof}
\begin{remark}\label{rk:Toen1}
Let $\cA$ be a spectral category. Since $[\cA\op]=[\cA]\op$, the opposite $F\op$ of a stable quasi-equivalence $F:\cA \to \cB$ is a stable quasi-equivalence. Therefore, if in Proposition~\ref{prop:Toen1}, we take $\cM= \Sp_{\bbU}$, we obtain an equivalence of triangulated categories $\bbL F_!:\cD(\cA\op) \stackrel{\sim}{\to} \cD(\cB\op)$, which restricts to an equivalence on its subcategories of compact objects $\bbL F_!:\cD_c(\cA\op) \stackrel{\sim}{\to} \cD_c(\cB\op)$. See \cite[4.2.7]{Neeman} for the notion of compact object.
\end{remark}
\begin{proposition}{(\cite[Prop.\,3.3]{Toen})}\label{prop:Toen2}
Let $\cA$ be a $\bbU$-small spectral category, such that for all objects $x,y \in \cA$ the symmetric spectrum $\cA(x,y)$ is cofibrant, and $\cM$ a $\bbV$-small cofibrantly generated $\Sp_{\bbU}$-model category. Then for any $x \in \cA$, the evaluation functor
$$
\begin{array}{rcl}
x^{\ast}: \cM^{\cA} & \too & \cM \\
M & \mapsto & M(x)
\end{array}
$$
preserves fibrations, cofibrations and weak equivalences.
\end{proposition}
\begin{proof}
The proof is analogous to the one of \cite[Prop.\,3.3]{Toen}. Replace the tensor product $-\otimes-$ by the smash product $-\wedge-$ and the notion of $C(k)$-model category (see \cite[\S 3]{Toen}) by the notion of $\Sp_{\bbU}$-model category (\ref{def:SpeMod}).
\end{proof}
\begin{remark}\label{rk:Toen2}
Let $\cB$ be a cofibrant spectral category. By Proposition~\ref{prop:cof}, for all objects $x, y \in \cB$, the symmetric spectrum $\cB(x,y)$ is cofibrant. Moreover, since the domains of the generating cofibrations in $\Sp_{\bbU}$ are cofibrant, the same holds for the generating cofibrations in $\Bmod_{\bbU}$. By Proposition~\ref{prop:Toen2}, if $M$ is cofibrant in $\Bmod_{\bbU}$ and $x \in \cB$, the symmetric spectrum $M(x)$ is cofibrant. This implies that if $Z \to Z'$ is a stable equivalence in $\Sp_{\bbU}$, so is $Z \wedge M(x) \to Z' \wedge M(x)$ (\cite[II-5.1]{Schwede}) and so the induced map $Z \wedge M \to Z' \wedge M$ is a weak equivalence in $\Bmod_{\bbU}$. In conclusion, the $\bbV$-small $\Sp_{\bbU}$-model category $\Bmod_{\bbU}$ (or $\Bomod_{\bbU}$) satisfies the conditions of Proposition~\ref{prop:Toen1}.
\end{remark}
Let $\cB$ be a $\bbU$-small fibrant spectral category. For every $x \in \cB$, the object $\hm_x \in \Bomod_{\bbU}$ is fibrant and cofibrant. Therefore we have a morphism in $\Spcat_{\bbV}$ (on the left)
$$
\begin{array}{lccr}
\hm_{-}: \cB \too \cBh & & \hm_{-}: [\cB] \too \cD(\cBo)\,,
\end{array}
$$
which induces the $\bbV$-small functor (on the right).
\begin{definition}\label{def:quasi-rep}
Let $\cB$ be a fibrant spectral category. A $\cBo$-module is called {\em quasi-representable} if it belongs to the essential image of the fully faithful functor
$$\underline{h}_{-}:[\cB] \too \cD(\cBo)\,.$$
\end{definition}
\subsection{Bimodules}
Given spectral categories $\cA$ and $\cB$, its {\em smash product} $\cA \wedge \cB$ is defined as follows: the set of objects of $\cA \wedge \cB$ is $\mbox{obj}(\cA)\times \mbox{obj}(\cB)$ and for two objects $(x,y)$ and $(x',y')$ in $\cA \wedge \cB$, we define
$$ (\cA \wedge \cB)((x,y),(x',y'))= \cA(x,x') \wedge \cB(y,y')\,.$$
This defines a symmetric monoidal structure on $\Spcat$, which is easily seen to be closed. However, the model structure of Theorem~\ref{thm:Model} endowed with the symmetric monoidal structure $-\wedge-$ is {\em not} a symmetric monoidal model category, as the smash product of two cofibrant objects in $\Spcat$ is not cofibrant in general. Nevertheless, the bi-functor $-\wedge-$ can be derived into a bi-functor
$$ - \Dsma - : \Ho(\Spcat) \times \Ho(\Spcat) \too \Ho(\Spcat)$$
defined by $\cA \Dsma \cB := Q(\cA) \wedge \cB$, where $Q(A)$ is a cofibrant replacement functor in $\Spcat$, which acts by the identity on the sets of objects (\ref{rk:cof}). 

Now, let $\cA$ and $\cB$ be spectral categories, with $\cB$ fibrant. For every object $x \in \cA$, there exists a spectral functor $\cBo \to \cA \wedge \cBo$ sending $y \in \cB$ to $(x,y)$, and
$$ \cBo(y,z) \too (\cA \wedge \cBo)((x,y),(x,z))=\cA(x,x) \wedge \cBo(y,z)$$
being the smash product of the unit $\bbS \to \cA(x,x)$ with the identity on $\cBo(y,z)$. As $\cA$ and $Q(\cA)$ have the same set of objects, one sees that for any $x \in \cA$, we have a natural spectral functor
$$ i_x: \cBo \too Q(\cA)\wedge\cBo = \cA \Dsma \cBo\,.$$
\begin{definition}\label{def:rep}
Let $\cA$ and $\cB$ be spectral categories, with $\cB$ fibrant. We denote by $\rep(\cA, \cB)$ the full subcategory of $\cD(\cA \Dsma \cBo)$, whose objects are the $\cA \Dsma \cBo$-modules $M$ such that $i_x^{\ast}(M)$ is quasi-representable (\ref{def:quasi-rep}), for all objects $x \in \cA$. We denote by $\mathsf{Iso}\,\rep(\cA, \cB)$ the set of isomorphism classes of $\rep(\cA,\cB)$.
\end{definition}
%-------------------------------------------------
\section{Homotopy category}\label{sec:chapter1}
%-------------------------------------------------
In this Chapter, we describe the Hom-sets in the homotopy category $\Ho(\Spcat)$, in terms of isomorphism classes of a certain derived category of bimodules, see Theorem~\ref{thm:Hom}.
\subsection{Path object}
\begin{definition}
Let $\cB$ be a $\bbU$-small fibrant spectral category. The $\bbV$-small spectral category $P(\cBh)$ (see~\ref{not:hat}) is defined as follows:
\begin{itemize}
\item[-] its objects are the fibrations $\phi:x \to y\times z$ in $\Bomod_{\bbU}$ (with $x,y,z \in \cBh$), whose components $x \stackrel{\sim}{\to} y$ and $x \stackrel{\sim}{\to} z$ are weak equivalences.
\item[-] given two objects $\phi: x \to y \times z$ and $\phi': x' \to y' \times z'$ in $P(\cBh)$, the symmetric spectrum $P(\cBh)(\phi, \phi')$ is defined by the following pullback square
$$
\xymatrix{
P(\cBh)(\phi, \phi') \ar[d] \ar[r] \ar@{}[dr]|{\ulcorner} & \cBh(x,x') \ar[d]^{\phi'_{\ast}} \\
\cBh(y,y') \times \cBh(z,z') \ar[r] & \underline{\Bomod_{\bbU}}(x, y'\times z')\,.
}
$$
\end{itemize}
\end{definition}
Since $x$ is cofibrant and $\phi'$ is a fibration in $\Bomod_{\bbU}$, the compatibility axiom (SP) implies that the morphism $\phi'_{\ast}$ is a stable projective fibration in $\Sp_{\bbU}$. Therefore, since fibrations are stable under base-change, we conclude that $P(\cBh)$ is a $\bbV$-small fibrant spectral category. We have two natural projection morphisms $\pi, \pi': P(\cBh) \to~\cBh$ in $\Spcat_{\bbV}$, given by 
$$ \pi(\phi:x \to y\times z) =y \,\,\,\,\,\,\,\,\,\, \pi'(\phi:x \to y\times z) =z\,.$$
By remark~\ref{rk:chunck}, we have also a morphism $\tau: \cBh \to P(\cBh)$ in $\Spcat_{\bbV}$, which maps an object $x \in \cBh$ to the map $q$ appearing in a chosen functorial factorization
$$ x \stackrel{p}{\too} \overline{x} \stackrel{q}{\too} x \times x$$
of the diagonal map, where $p$ is a trivial cofibration and $q$ is a fibration. We obtain in this way a commutative diagram in $\Spcat_{\bbV}$
$$
\xymatrix{
\cBh \ar[rr]^{\Delta} \ar[dr]_{\tau} && \cBh \times \cBh \\
& P(\cBh) \ar[ur]_{\pi\times \pi'} & \,.
}
$$
\begin{theorem}\label{thm:path}
The $\bbV$-small spectral category $P(\cBh)$ is a path object~\cite[7.3.2(3)]{Hirschhorn} for $\cBh$, with respect to the model structure of Theorem~\ref{thm:Model}.
\end{theorem}
\begin{proof}
We start by showing that $\tau$ is a stable quasi-equivalence. By the two-out-of-three property, it is enough to show that $\pi:P(\cBh) \to \cBh$ is a stable quasi-equivalence. Since $\tau$ is a section of $\pi$, $\pi$ satisfies condition S2'). We now show that it also satisfies condition S1), \ie for all objects $\phi: x\to y \times z$ and $\phi':x' \to y'\times z'$ in $P(\cBh)$, the induced morphism 
$$ P(\cBh)(\phi, \phi') \too \cBh(y,y')$$
is a stable equivalence. Notice that we have a commutative diagram
$$
\xymatrix{
P(\cBh)(\phi, \phi') \ar[d] \ar[r] \ar@{}[dr]|{\ulcorner} & \cBh(x,x') \ar[d]^{\phi'_{\ast}} \\
\cBh(y,y') \times \cBh(z,z') \ar[r] \ar[d]_{\sim} & \underline{\Bomod_{\bbU}}(x, y'\times z') \ar[d]^{\sim} \\
\cBh(y,y') \times \cBh(x,z') \ar[r] \ar[d] \ar@{}[dr]|{\ulcorner} & \cBh(x,y') \times \cBh(x,z') \ar[d] \\
\cBh(y,y') \ar[r] & \cBh(x,y') \,.
}
$$
Since $x \stackrel{\sim}{\to} y$ is a weak equivalence between cofibrant objects in $\Bomod_{\bbU}$ and $y'$ is fibrant, the bottom horizontal morphism is a stable equivalence in $\Sp_{\bbU}$.  Moreover, since all the above squares are homotopy cartesian, the upper horizontal morphism is also a stable equivalence in $\Sp_{\bbU}$. Finally, since in the commutative square
$$
\xymatrix{
P(\cBh)(\phi, \phi') \ar[d] \ar[rr]^{\sim} && \cBh(x,x') \ar[d] \\
\cBh(y,y') \ar[rr]_{\sim} && \cBh(x,y')
}
$$
the right vertical morphism is a stable equivalence in $\Sp_{\bbU}$, the spectral functor $\pi$ satisfies condition S1) and so $\tau$ is a stable quasi-equivalence. It remains to prove that
$$ \pi\times \pi': P(\cBh) \too \cBh$$
is a fibration. Since $P(\cBh)$ and $\cBh$ are $\bbV$-small fibrant spectral categories, it is enough to show that $\pi\times \pi'$ is a levelwise fibration, see \cite[Prop.\,4.15]{Spectral}. Condition F1) follows from the fact that $\phi'_{\ast}$ is a stable projective fibration in $\Sp_{\bbU}$. In what concerns condition F2), we need to show that the $\bbV$-small simplicial functor
$$ P(\cBh)_0 \too \cBh_0 \times \cBh_0=(\cBh \times \cBh)_0$$
is a fibration. Notice that $\cBh_0$ is the simplicial category of fibrant and cofibrant objects in $\Bomod_{\bbU}$. Therefore, given an object $\phi:x \to y \times z$ in $P(\cBh)_0$ and weak equivalences $f:y \stackrel{\sim}{\to} y'$ and $g:z \stackrel{\sim}{\to} z'$ in $\cBh_0$, we can factor the composite map $x \to y' \times z'$ as a trivial cofibration followed by a fibration $\phi': x' \to y' \times z'$. We obtain in this way a morphism $\alpha: \phi \to \phi'$ in $P(\cBh)_0$. Since $\pi_0:P(\cBh)_0 \to \cBh_0$ is a DK-equivalence and $\pi_0(\alpha)=f$, we conclude that $\alpha$ becomes invertible in homotopy category of $P(\cBh)_0$. 
\end{proof}
\subsection{Hom-sets and bimodules}
\begin{proposition}\label{prop:Lurie}
Let $\cA$ and $\cB$ be a $\bbU$-small spectral categories, with $\cA$ cofibrant and $\cB$ fibrant. Then the natural map
$$ 
\begin{array}{lcccr}
[\,\cA, \cBh\,]  \too  \mbox{Iso}\,\cD(\cA \wedge \cBo) &&&
F  \mapsto  [\,(x,y) \mapsto \underline{\Bomod_{\bbU}}(y,F(x))\,]\\
\end{array}
$$
is well-defined and injective.
\end{proposition}
\begin{proof}
Suppose first that $F$ and $F'$ coincide in $[\cA, \cBh]$. Then by Theorem~\ref{thm:path}, we have a commutative diagram
$$
\xymatrix{
& \cBh \\
\cA \ar[ur]^F \ar[r]^H \ar[dr]_{F'} & P(\cBh) \ar[u]_{\pi} \ar[d]^{\pi'} \\
& \cBh \,. 
}
$$
Notice that the homotopy $H$ furnishes us a new morphism $F'':\cA \to \cBh$ in $\Spcat_{\bbV}$, equiped with weak equivalences $F'' \to F$ and $F'' \to F'$. This implies that the natural map is well defined.

Now, suppose that $F$ and $F'$ coincide in $\mbox{Iso}\,\cD(\cA\wedge \cBo)$. Let $\alpha:F''\to F$ be a cofibrant resolution of $F$. Since $\cA$ is cofibrant and $F''$ is fibrant and cofibrant in $(\Bomod_{\bbU})^{\cA}$, Propositions~\ref{prop:cof} and~\ref{prop:Toen2} imply that $F''$ take values in $\cBh$ and so it corresponds to a morphism $F'':\cA \to \cBh$ in $\Spcat_{\bbV}$. Since $F$ and $F'$ are weakly equivalent and $F''$ is cofibrant, there exists a weak equivalence $\beta: F'' \to F'$. Finally, we can factor the map $\alpha \times \beta$ as $F'' \stackrel{u}{\to} F''' \stackrel{v}{\to} F \times F'$, where $u$ is a trivial cofibration and $v$ a fibration. The map $v$ can be viewed as an object of $P(\cBh)$ and so it gives rise to an homotopy from $F$ to $F'$. 
\end{proof}
\begin{theorem}\label{thm:Hom}
Given $\bbU$-small spectral categories $\cA$ and $\cB$, with $\cB$ fibrant, we have a natural bijection
$$[\cA, \cB] \stackrel{\sim}{\too} \mathsf{Iso}\,\rep(\cA, \cB)\,.$$
\end{theorem}
\begin{proof}
We can assume that $\cA$ is cofibrant and so that $Q(\cA)=\cA$. Since the Yoneda morphism $\underline{h}_-: \cB \to \cBh$ in $\Spcat_{\bbV}$ is fully faithfull it induces, as in \cite[Corollary\,2.4]{Toen}, an injective map $[\cA, \cB] \to [\cA, \cBh]$. By composing it with the one of Proposition~\ref{prop:Lurie}, we obtain an injective map$$[\cA, \cB] \too  \mbox{Iso}\,\cD(\cA\wedge \cBo)\,,$$
which factors through $\mathsf{Iso}\,\rep(\cA,\cB)$. Now let $F$ be an object of $\mathsf{Iso}\,\rep(\cA, \cB)$, which we can assume to be fibrant and cofibrant. Since $\cA$ is cofibrant, $F$ corresponds by Proposition~\ref{prop:Toen2}, to a morphism
$$ F: \cA \too \cBh^{qr}$$
in $\Spcat_{\bbV}$, where $\cBh^{qr}\subset \cBh$ denotes the full spectral subcategory of quasi-representable $\bbU$-small $\cBo$-modules (\ref{def:quasi-rep}). Notice that we have a diagram in $\Spcat_{\bbV}$
$$
\xymatrix{
\cA \ar[r]^F & \cBh^{qr} \\
& \cB \ar[u]_{\underline{h}_{-}} \,,
}
$$
where $\underline{h}_{-}$ a stable quasi-equivalence. Since $\cA$ is cofibrant and $\cB$ fibrant, there exists a morphism $G:\cA \to \cB$ in $\Spcat_{\bbU}$ such that $F$ and $\underline{h}_{-} \circ G$ are homotopic. By Proposition~\ref{prop:Lurie}, we conclude that they are therefore isomorphic in $\rep(\cA, \cB)$ and so the proof is finished.
\end{proof}
%-----------------------------------------------------------------------
\section{Triangulated envelope}\label{sec:chapter3}
%-----------------------------------------------------------------------
\subsection{Triangulated spectral categories}\label{sec:chapter2} By inspiring ourselves in To{\"e}n's lectures~\cite{Toen1} on dg categories, we introduce the following notion.
\begin{definition}\label{def:triang}
A $\bbU$-small fibrant spectral category $\cA$ is {\em triangulated} if the $\bbV$-small functor
$$\hm_-: [\cA] \stackrel{\sim}{\too} \cD_c(\cAo)$$
is an equivalence.
\end{definition}
\begin{notation}
We denote by $\Ho(\Spcat)^{tr} \subset \Ho(\Spcat)$ the full subcategory of triangulated spectral categories.
\end{notation}
\begin{remark}\label{rem:point}
If $\cA \in \Ho(\Spcat)^{tr}$, then $[\cA]$ is a (idempotent complete~\cite{Balmer}) triangulated category and if $f:\cA \to \cB$ is a morphism in $\Ho(\Spcat)^{tr}$, the induced functor $[f]: [\cA] \to [\cB]$ is triangulated. Moreover, we can always represent an object $\cA$ in $\Ho(\Spcat)^{tr}$ by a triangulated spectral category endowed with a zero object: consider the spectral category $\cA_+$ obtained from $\cA$ by adding a {\em zero object} $0$, \ie $\cA_+(x,0)=\ast$ and $\cA_+(0,x)=\ast$ for every $x \in \cA$. We have a natural fully faithful spectral functor $\cA \to \cA_+$ and under the equivalence $\hm_-:[\cA] \stackrel{\sim}{\to} \cD_c(\cA\op)$, the trivial $\cA\op$-module corresponds to an object of $[\cA]$, which becomes isomorphic in $[\cA_+]$ to the zero object. Therefore $\cA$ and $\cA_+$ are isomorphic in $\Ho(\Spcat)^{tr}$. 
\end{remark}
\begin{theorem}\label{thm:envelope}
The natural inclusion functor $\Ho(\Spcat)^{tr} \subset \Ho(\Spcat)$ admits a left adjoint $(-)^{\sharp}_{pe}$, which we refer to as the {\em triangulated envelope}.
\end{theorem}
According to \cite[Thm.\,2(ii)-IV]{MacLane}, to prove Theorem~\ref{thm:envelope} it is enough to construct for every $\cA \in \Ho(\Spcat)$ a morphism $ \theta: \cA \too \Ape$ with $\Ape \in \Ho(\Spcat)^{tr}$, such that for every $\cB \in \Ho(\Spcat)^{tr}$, the induced map $$\theta^{\ast}: [\Ape, \cB] \stackrel{\sim}{\too} [\cA, \cB]$$ is bijective.
Theorem~\ref{thm:envelope} will follow from Propositions~\ref{prop:envelope} and \ref{prop:envelope1} below.

\subsection{Construction of the triangulated envelope} Let $\cA \in \Ho(\Spcat_{\bbU})$, which we can assume to be a $\bbU$-small fibrant and cofibrant spectral category. Consider the fully faithful Yoneda morphism in $\Spcat_{\bbV}$:
$$ \hm_-: \cA \too \cAh \subset \Aomod_{\bbU} \,.$$
Let $\cA_{pe}$ be the $\bbV$-small full spectral subcategory of $\cAh$, whose objects are the $\bbU$-small {\em perfect} $\cAo$-modules, \ie those $\cAo$-modules that become compact in $\cD(\cAo)$. Notice that we have an equivalence $[\cA_{pe}] \simeq \cD_c(\cAo)$ and moreover every object $z \in \cA_{pe} \subset \Aomod_{\bbU}$ is equivalent in $\cD_c(\cAo)$ to the retract of some $\bbU$-small $\cAo$-module $z'$, obtained by a finite composition
$$ z_0=\ast \to z_1 \to \cdots \to z_i \to z_{i+1} \to \cdots \to z_r=z'\,,$$
such that for every $i$, we have a (homotopy) pushout diagram
$$
\xymatrix{
F_n\partial\Delta[m]^+ \Dsma \hm_{x_i} \ar[r] \ar[d] \ar@{}[dr]|{\lrcorner} & z_i \ar[d] \\
F_n\Delta[m]^+ \ar[r] \Dsma \hm_{x_i} & z_{i+1}\,,
}
$$ 
where $x_i \in \cA$ and $F_n\partial \Delta[m]^+ \to F_n \Delta[m]^+$ is a generating cofibration in $\Sp_{\bbU}$. Therefore, we can take inside $\cA_{pe}$ the smallest $\bbU$-small full spectral subcategory $\Ape \subset \cA_{pe}$, which contains a representative of each isomorphism class in $\cD_c(\cAo)$. Notice that we also have an equivalence $[\Ape] \simeq \cD_c(\cAo)$. In conclusion, we have constructed a fully-faithful morphism in $\Spcat_{\bbU}$ 
$$ \cA \too \Ape \subset \cA_{pe} \subset \cAh \subset \Aomod_{\bbU} \,.$$
\begin{notation}\label{not:envelope}
We denote by 
$$ \theta: \cA \too \Ape \,.$$
the corresponding morphism in $\Ho(\Spcat_{\bbU})$. The $\bbU$-small spectral category $\Ape$ is called the {\em triangulated envelope of $\cA$}.
\end{notation}
\begin{lemma}{(\cite[Lemma\,7.5]{Toen})}\label{lem:general}
Let $\cM$ be a $\bbV$-small cofibrantly generated $\Sp_{\bbU}$-model category (\ref{def:SpeMod}), which satisfies the conditions of Proposition~\ref{prop:Toen1}. Then:
\begin{itemize}
\item[(1)] The Quillen adjunction
$$
\begin{array}{lccr}
\theta_!: \cM^{\cA} \to \cM^{\Ape} && \cM \leftarrow \cM^{\Ape}: \theta^{\ast}
\end{array}
$$
is a Quillen equivalence.
\item[(2)] For any $F \in \cM^{\Ape}$, and any $\bbU$-small diagram $X: I \to \Aomod_{\bbU}$ of perfect and cofibrant objects in $\Aomod_{\bbU}$, the natural morphism
$$ \mbox{hocolim}_i \,F(X_i) \stackrel{\sim}{\too} F(\mbox{hocolim}_i\,X_i)$$
is an isomorphism in $\Ho(\cM)$.
\item[(3)] For any $F \in \cM^{\Ape}$, $x \in \cA$ and $Z \in \{ F_n\partial \Delta[m]^+, F_n\Delta[m]^+\}_{n,m \geq 0}$, the natural morphism
$$ Z \Dsma F(\hm_x) \stackrel{\sim}{\too} F(Z \Dsma \hm_x)$$
is an isomorphism in $\Ho(\cM)$. 
\end{itemize}
\end{lemma}
\begin{proof}
The proof is analogous to the one of \cite[Lemma\,7.5]{Toen}. Replace the tensor product $-\otimes-$ by the smash product $-\wedge-$ and the notion of $C(k)$-model category (see \cite[\S 3]{Toen}) by the notion of $\Sp_{\bbU}$-model category (\ref{def:SpeMod}).
\end{proof}
\subsection{Proof of Theorem~\ref{thm:envelope}}
\begin{proposition}\label{prop:envelope}
Let $\cB \in \Ho(\Spcat_{\bbU})^{tr}$. Then the induced map
$$ \theta^{\ast}: [\Ape, \cB] \stackrel{\sim}{\too} [\cA, \cB]$$
is bijective.
\end{proposition}
\begin{proof}
We can assume that $\cB$ is a $\bbU$-small fibrant and cofibrant spectral category. By Remark~\ref{rk:Toen2}, the $\Sp_{\bbU}$-model category $\cM= \Bomod_{\bbU}$ satisfies the conditions of Proposition~\ref{prop:Toen1} and so of Lemma~\ref{lem:general}. We obtain then a (derived) equivalence
$$
\xymatrix{
\Ho((\Bomod_{\bbU})^{\Ape}) \ar[d]_{\sim} \ar[rr]^{\theta^{\ast}}_{\sim} && \Ho((\Bomod_{\bbU})^{\cA}) \ar[d]^{\sim}\\
\cD((\Ape) \Dsma \cB\op) \ar[rr]_{\sim} && \cD(\cA \Dsma \cB\op) \,.
}
$$
Now, consider the following commutative diagram (\ref{def:rep})
$$ 
\xymatrix{
\rep(\Ape, \cB) \ar[d] \ar@{^{(}->}[rr] && \cD(\Ape \Dsma \cBo) \ar[d]^{\sim} \\
\rep(\cA, \cB) \ar@{^{(}->}[rr] && \cD(\cA \Dsma \cBo)\,.
}
$$
Since the right vertical functor is an equivalence, the left vertical functor is fully faithful. We now show that it is also essentially surjective. Let $F \in \rep(\cA, \cB)$. Its image in $\cD(\cA \Dsma \cBo)$ comes an object $\overline{F} \in \Ho((\cBo\text{-}\mbox{Mod}_{\bbU})^{\Ape})$. For every $x_i \in \cA$, the $\cBo$-module $\overline{F}(\hm_{x_i})=F(x_i)$ belongs to $\cD_c(\cBo)$ and by Lemma~\ref{lem:general} so it does $\overline{F}(Z_i \Dsma \hm_{x_i})$, with $Z_i \in \{F_n\partial\Delta[m]^+, F_n\Delta[m]^+ \}_{n,m \geq 0}$. Since every object $z \in \Ape$ can be constructed as a finite (homotopy) colimit of objects of the form $Z_i \Dsma \hm_{x_i}$, with $x_i \in \cA$ and $Z_i \in \{F_n\partial\Delta[m]^+, F_n\Delta[m]^+ \}_{n,m \geq 0}$, we conclude that $\overline{F}(z)$ belongs also to $\cD_c(\cBo)$. Since $\cB \in \Ho(\Spcat_{\bbU})^{tr}$, we have an equivalence of categories $[\cB] \simeq \cD_c(\cBo)$, which implies that $\overline{F} \in \rep(\Ape, \cB)$. This shows that the left vertical functor is essentially surjective. Finally, by Theorem~\ref{thm:Hom}, we conclude that the induced map
$$\theta^{\ast}: [\Ape, \cB] \stackrel{\sim}{\too} [\cA, \cB]\,.$$
is bijective.
\end{proof}
\begin{remark}\label{rk:eqrep}
The proof of Proposition~\ref{prop:envelope} shows us that the bijection (on the left)
$$
\begin{array}{lccr}
[\Ape, \cB] \stackrel{\sim}{\too} [\cA, \cB] & && \rep(\Ape, \cB) \stackrel{\sim}{\too} \rep(\cA, \cB)\,,
\end{array}$$
follows from the equivalence (on the right).
\end{remark}
\begin{proposition}\label{prop:envelope1}
Let $\cA$ be an object in $\Ho(\Spcat_{\bbU})$. Then its triangulated envelope $\Ape$ belongs to $\Ho(\Spcat_{\bbU})^{tr}$.
\end{proposition}
\begin{proof}
We have a fully faithful functor
$$ \hm_{-}: [\Ape] \too \cD_c((\Ape)\op) \,,$$
which by Lemma~\ref{lem:general} is moreover triangulated. Recall that $\cD_c((\Ape)\op)$ is the smallest thick triangulated subcategory of $\cD((\Ape)\op)$, which contains the images of $\hm_{-}$. Since by construction the category $[\Ape]$ is stable under (co)suspensions, homotopy colimits and is moreover idempotent complete, we conclude that $\hm_{-}$ is also essentially surjective.
\end{proof}
%\subsection{Semi-additive structure}
%\begin{proposition}\label{prop:Qadd}
%The category $\Ho(\Spcat)^{tr}$ has products ($-\times-$) and co-products ($-\underset{pe}{\amalg}-$). Moreover, for all $\cA_1, \cA_2 \in \Ho(\Spcat)^{tr}$, the natural morphism
%$$ \cA_1 \underset{pe}{\amalg} \cA_2 \stackrel{\sim}{\too} \cA_1 \times \cA_2$$
%is an isomorphism.
%\end{proposition}
%\begin{remark}
%Proposition~\ref{prop:Qadd} implies that the set $[\cT_1, \cT_2]$ is endowed with a structure of commutative monoide, making the composition operation in $\Ho(\Spcat)^{tr}$ bilinear.
%\end{remark}
%\begin{proof}
%The homotopy category $\Ho(\Spcat)$ has products ($-\times-$) and co-products ($-\amalg-$). These are obtained by deriving the products and coproducts in $\Spcat$. Clearly $\Ho(\Spcat)^{tr} \subset \Ho(\Spcat)$ is stable under products. Now, let $\cA_1, \cA_2 \in \Ho(\Spcat)^{tr}$. By adjunction (\ref{thm:envelope}), the coproduct of $\cA_1$ and $\cA_2$ corresponds to
%$$ \cA_1 \underset{pe}{\amalg}\cA_2 := (\cA_1 \amalg \cA_2)^{\sharp}_{pe}\,.$$
%Finally, notice that we have the following isomorphisms in $\Ho(\Spcat)^{tr}$
%$$ \cA_1 \underset{pe}{\amalg}\cA_2 := (\cA_1 \amalg \cA_2)^{\sharp}_{pe} \simeq (\cA_1)_{pe}^{\sharp} \times (\cA_1)_{pe}^{\sharp} \simeq \cA_1 \times \cA_2\,.$$
%\end{proof}

%-----------------------------------------------------------------------
\section{Morita equivalences}\label{sec:chapter4}
%-----------------------------------------------------------------------
In this Chapter we describe the localization of $\Spcat$ with respect to the following class of spectral functors:
\begin{definition}\label{def:Morita}
A spectral functor $F:\cA \to \cB$ is called a {\em Morita equivalence} if its derived extension of scalars $$ \bbL F_!: \cD(\cAo) \stackrel{\sim}{\too} \cD(\cBo)$$
is an equivalence.  A morphism $f:\cA \to \cB$ in $\Ho(\Spcat)$ is called a {\em Morita morphism} if $f^{\sharp}_{pe}: \Ape \to \Bpe$ is an isomorphism in $\Ho(\Spcat)^{tr}$.
\end{definition}
\begin{remark}\label{rk:Morita}
Since the functor $\bbL F_!$ commutes with sums and the triangulated categories $\cD(\cAo)$ and $\cD(\cBo)$ are compactly generated~\cite[\S 8.1]{Neeman}, $F$ is a Morita equivalence if and only if the induced functor $\bbL F_!: \cD_c(\cAo) \to \cD_c(\cBo)$ is an equivalence.
\end{remark}
\begin{lemma}\label{lem:triang}
Let $f:\cA \to \cB$ be a morphism in $\Ho(\Spcat)^{tr}$. Then $f$ is a isomorphism if and only if the induced triangulated functor
$$ [f]: [\cA] \stackrel{\sim}{\too} [\cB]$$
is an equivalence.
\end{lemma}
\begin{proof}
If $f$ is an isomorphism, $[f]$ is clearly an equivalence. Let us now prove the converse. We can assume that both $\cA$ and $\cB$ are fibrant and cofibrant, and so we can represent $f$ by a spectral functor $F:\cA \to \cB$. Therefore, $f$ is an isomorphism in $\Ho(\Spcat)^{tr}$ if and only if it $F$ is a stable quasi-equivalence (\ref{def:quasi}). By hypothesis, $F$ satisfies condition S2). Since $\cA$ and $\cB$ are fibrant, to verify condition S1), it is enough by Proposition~\ref{prop:fibrant} to show that for all objects $x,y \in \cA$, the morphism in $\Sp$
$$ F(x,y): \cA(x,y) \too \cB(Fx,Fy)$$
induces an isomorphism in all stable homotopy groups $\pi^s_j(-), j \in \bbZ$. We have commutative squares
$$
\xymatrix{
\cA \ar[r]^F \ar[d]_{\hm_-} & \cB \ar[d]^{\hm_-} && [\cA] \ar[r]^{[F]}_{\sim} \ar[d]_{\sim} & [\cB] \ar[d]^{\sim} \\
\Aomod\ar[r]_{F_!} & \Bomod&& \cD_c(\cAo) \ar[r]_{\bbL F_!}^{\sim} & \cD_c(\cBo)\,.
}
$$
By~\cite[Lemma\,3.5.2]{SS} and the above squares, we have a commutative diagram
$$
\xymatrix{
\pi^s_j \cA(x,y) \ar[d]_{\sim} \ar[rr] && \pi^s_j \cB(Fx,Fy) \ar[d]^{\sim} \\
\pi^s_j \underline{\Aomod}(\hm_x, \hm_y) \ar[rr] \ar[d]_{\sim} && \pi^s_j \underline{\Bomod} (\hm_{Fx}, \hm_{Fy})\ar[d]^{\sim} \\
[\Sigma^j(\hm_x), \hm_y] \ar[rr]_{\sim} && [\Sigma^j(\hm_{Fx}), \hm_{Fy}]\,,
}
$$
where $\Sigma^j(-), j \in \bbZ$ denotes the $j$-th suspension functor. Since the lower horizontal map is an isomorphism, we conclude that condition S1) is verified and so that $F$ is a stable quasi-equivalence.
\end{proof}
\begin{proposition}\label{prop:Morita}
A spectral functor $F:\cA \to \cB$ in $\Spcat$ is a Morita equivalence if and only if its induced morphism $f:\cA \to \cB$ in $\Ho(\Spcat)$ is a Morita morphism.
\end{proposition}
\begin{proof}
By Lemma~\ref{lem:triang}, $f^{\sharp}_{pe}: \Ape \too \Bpe$ is an isomorphism in $\Ho(\Spcat)^{tr}$ if and only if the induced triangulated functor $[f^{\sharp}_{pe}]: [\Ape] \to [\Bpe]$ is an equivalence. Since we have the (up to equivalence) commutative square
$$
\xymatrix{
[\Ape] \ar[d]_{\sim} \ar[rr]^{[f_{pe}^{\sharp}]} && [\Bpe] \ar[d]^{\sim} \\
\cD_c(\cAo) \ar[rr]_{\bbL F_!} && \cD_c(\cBo)\,,
}
$$
Remark~\ref{rk:Morita} allows us to conclude the proof.
\end{proof}
%\begin{lemma}\label{lem:local}
%Let $\cC \in \Ho(\Spcat)$. Then $\cC \in \Ho(\Spcat)^{tr}$ if and only if for every Morita morphism $f: \cA \to \cB$ (Def.\,\ref{def:Mmorph}), the induced map
%$$f^{\ast}: [\cB, \cC] \stackrel{\sim}{\too} [\cA, \cC]$$
%is bijective.
%\end{lemma}
%\begin{proof}
%If $\cC \in \Ho(\Spcat)^{tr}$, proposition~\ref{prop:Morita} implies that the induced map is bijective. Let us now prove the converse. Consider the triangulated envelope $\eta_{\cC}: \cC \too \cC^{\sharp}_{pe}$ of $\cC$ (\ref{not:envelope}).
%Since the induced map
%$$ \theta^{\ast}: [\cC^{\sharp}_{pe}, \cC] \stackrel{\sim}{\too} [\cC, \cC]$$
%is bijective, $\cC$ is a retract of $\cC^{\sharp}_{pe}$ in $\Ho(\Spcat)$. By proposition~\ref{prop:envelope} $\Cpe \in \Ho(\Spcat)^{tr}$ and so by lemma~\ref{lem:factors} we conclude that $\cC \in \Ho(\Spcat)^{tr}$.
%\end{proof}
\begin{theorem}\label{thm:localization}
The composition
$$ \Spcat \too \Ho(\Spcat) \stackrel{(-)^{\sharp}_{pe}}{\too} \Ho(\Spcat)^{tr} $$
is the localization functor associated to the Morita equivalences.
\end{theorem}
\begin{proof}
Let $\Psi: \Spcat \to \mathsf{D}$ be a functor which sends Morita equivalences to isomorphisms. Since by Remark~\ref{rk:Toen1} every quasi-equivalence is a Morita equivalence, the functor $\Psi$ descends to the homotopy category $\Ho(\Spcat)$. It is a general fact (see~\cite[Prop.\,1.3]{Zisman}), that the left adjoint functor
$$(-)^{\sharp}_{pe}: \Ho(\Spcat) \too \Ho(\Spcat)^{tr}$$
is the localization of $\Ho(\Spcat)$ with respect to the Morita morphisms. Since every morphism in $\Ho(\Spcat)$ can be represented by a spectral functor, Proposition~\ref{prop:Morita} shows us that the class of Morita equivalences in $\Spcat$ and the class of Morita morphisms in $\Ho(\Spcat)$ correspond one to another, under the localization functor
$$ \Spcat \too \Ho(\Spcat)\,.$$
This finishes the proof. 
\end{proof}
\section{Upper triangular matrices}\label{se:add}
%-----------------------------------------------------------------------

\begin{definition}\label{UTM}
An {\em upper triangular matrix} $\underline{M}$ is given by
$$
\begin{array}{rcl}
\underline{M} & := & \begin{pmatrix} \cA & X \\ \ast &
  \cC \end{pmatrix}\,,
\end{array}
$$
where $\cA$ and $\cC$ are spectral categories and $X$
is a $\cA\wedge \cCo$-module. The {\em totalization} $|\underline{M}|$ of $\underline{M}$ is the spectral category whose set of objects is the disjoint union of the sets of objects of $\cA$
and $\cC$ and whose morphisms are given by
$$
\begin{array}{rcl}
|\underline{M}|(x,y) & := & 
\left\{
\begin{array}{ccl}
\cA(x,y) & \mbox{if} & x,\, y \in \cA \\
\cC(x,y) & \mbox{if} & x,\, y \in \cC \\
X(x,y) & \mbox{if} & x \in \cA,\, y \in \cC \\
\ast & \mbox{if} &  x \in \cC,\, y \in \cA
\end{array}\right.
\end{array}\,.
$$
The composition is induced by the composition on $\cA$, $\cC$ and the $\cA\text{-}\cC$-bi-module structure on $X$. We have two natural inclusion spectral functors
$$
\begin{array}{rcl}
i_1: \cA \too |\underline{M}| & & i_2: \cC \too |\underline{M}| \,.
\end{array}
$$
\end{definition}
Let $\cI$ be the spectral category 
$$\xymatrix{
    1 \ar@(ul,dl)[]_{\bbS} \ar@/^/[r]^{\bbS}  &
    2 \ar@(ur,dr)[]^{\bbS} \ar@/^/[l]^{\ast} 
    }
$$
with two objects $1$ and $2$ such that $\cI(1,1)=\bbS$, $\cI(2,2)=\bbS$, $\cI(1,2)=\bbS$, $\cI(2,1)=\ast$, and composition given by multiplication. Given a spectral category $\cA$, we note by $T(\cA)$ the spectral category $\cA \wedge \cI$. Notice that $T(\cA)$ corresponds to the totalization of the upper triangular matrix
$$
\begin{pmatrix} \cA & \cA(-,-) \\ \ast &
  \cA
\end{pmatrix}\,.
$$
We have two natural inclusions
$$
\begin{array}{lccr}
\begin{array}{rcl}
i_1: \cA & \too & T(\cA) \\
x& \mapsto & (x,1)
\end{array}
& &
\begin{array}{rcl}
i_2: \cA & \too & T(\cA) \\
x & \mapsto & (x,2)
\end{array}
\end{array}
$$
and a projection spectral functor
$$P: T(\cA) \too \cA \,,$$
verifying the identities $P\circ i_1 = \mbox{Id}$ and $P\circ i_2 = \mbox{Id}$. Since a spectral functor from $\cI$ to a spectral category $\cB$, corresponds to specifying two objects $x$ and $y$ in $\cB$ plus a $0$-simplex in the degree zero component of the symmetric spectrum $\cB(x,y)$, the category 
$ \TAomod$ identifies with the category of morphisms in $\Aomod$. Therefore, we 
obtain the following extension of scalars functors:
$$
\begin{array}{lccr}
\xymatrix{
\TAomod  \ar[d]_{P_!} \\
\Aomod 
} & & & 
\xymatrix{
\TAomod \\
\Aomod \ar@<-1ex>[u]_{{i_2}_!} \ar@<1ex>[u]^{{i_1}_!} \,, 
}
\end{array}
$$
with
$$
\begin{array}{lcr}
{i_1}_!: L \mapsto (\ast \to L) && {i_2}_!: L \mapsto (L =L)
\end{array}
$$
and
$$
P_!:(X \to Y)\mapsto Y\,.
$$
\begin{remark}\label{rk:adjoints1}
Let $\cB$ be an object of $\Ho(\Spcat)^{tr}$. By Remark~\ref{rem:point}, we can assume that $\cB$ is a triangulated spectral category endowed with a zero object $0$. Therefore we have two new projection spectral functors:
$$
\begin{array}{lccr}
\begin{array}{rcl}
P_1: T(\cB) & \too & \cB \\
(x,1)& \mapsto & x\\
(x,2) & \mapsto & 0
\end{array}
& &
\begin{array}{rcl}
P_2: T(\cB) & \too & \cB \\
(x,1) & \mapsto & 0\\
(x,2) & \mapsto & x\,,

\end{array}
\end{array}
$$
which induce the following extensions of scalars
$$
\begin{array}{lccr}
\begin{array}{rcl}
{P_1}_! : \TBomod & \too & \Bomod \\
(X \to Y) & \mapsto & Y/X
\end{array}
& &
\begin{array}{rcl}
{P_2}_! : \TBomod & \too & \Bomod \\
(X \to Y) & \mapsto & X
\end{array} \,.
\end{array}
$$
In conclusion, if $\cB \in \Ho(\Spcat)^{tr}$, we have the equalities
$$
\begin{array}{cccc}
P_1 \circ i_1= \mbox{Id} & P_1 \circ i_2=\ast & P_2 \circ i_1=\ast & P_2\circ i_2=\mbox{Id}
\end{array}
$$
and
$$
\begin{array}{cc}
P \circ i_1=\mbox{Id} & P \circ i_2=\mbox{Id}\,.
\end{array}
$$
\end{remark}
\section{Universal matrix invariant}\label{ssec:additivization}
%-----------------------------------------------------------------------
\begin{definition}\label{def:additive}
Let $F: \Spcat \too \A$ be a functor with values in an additive category $\A$. We say that $F$ is a {\em matrix invariant of spectral categories} if it verifies the following two conditions:
\begin{itemize}
\item[M)] the functor $F$ sends the Morita equivalences (\ref{def:Morita}) to isomorphisms in $\A$;
\item[MA)] for every upper triangular matrix $\underline{M}$ (\ref{UTM}), the inclusion spectral functors induce an isomorphism in $\cA$
$$ [F(i_1) \,\, F(i_2)]: F(\cA) \oplus F(\cC) \stackrel{\sim}{\too} F(|\underline{M}|)\,.$$
\end{itemize}
\end{definition}
In this Chapter, we will construct the universal matrix invariant of spectral categories, see Theorem~\ref{thm:main}. \subsection{Additivization}
Notice first that given spectral categories $\cA$ and $\cB$, with $\cB$ a triangulated spectral category (\ref{def:triang}), the category $\rep(\cA,\cB)$ (\ref{def:rep}) is naturally triangulated. Moreover, given fibrant spectral categories $\cA$, $\cB$ and $\cC$, the bijection of Theorem~\ref{thm:Hom}, shows us that the composition in $\Ho(\Spcat)^{tr}$ is induced by the derived triangulated bi-functor
$$ 
\begin{array}{ccc}
- \Dsma_{\cB} - : \rep(\cA, \cB) \times \rep(\cB, \cC) & \too & \rep(\cA, \cC) \\
(X,Y) & \mapsto & X \Dsma_{\cB} Y\,.
\end{array}
$$
\begin{definition}\label{def:add}
Let $\Add$ be the category of fibrant objects in $\Ho(\Spcat)^{tr}$, whose Hom-sets are given by
$$ \Add(\cA, \cB):= K_0\rep(\cA,\cB)\,,$$
where $K_0\rep(\cA, \cB)$ denotes the Grothendieck group of $\rep(\cA, \cB)$. The composition is the induce one.
\end{definition}
We have a natural functor
$$ \Ho(\Spcat)^{tr} \too \Add$$
which sends each isomorphism class of $\rep(\cA, \cB)$ to the corresponding class in the Grothendieck group $K_0\rep(\cA,\cB)$.
\begin{lemma}
The category $\Add$ is additive.
\end{lemma}
\begin{proof}
By construction, the Hom-sets in $\Add$ are abelian groups and the composition operation is bilinear. Hence, it is enough by \cite[Thm.\,2-VII]{MacLane} to show that $\Add$ has direct sums. Given $\cA, \cB, \cC \in \Ho(\Spcat)^{tr}$, we have equivalences (see~\ref{rk:eqrep})
$$ \rep((\cA \amalg \cB)^{\sharp}_{pe}, \cC) \simeq \rep(\cA \amalg \cB, \cC) \simeq  \rep(\cA,\cC)\times \rep(\cB,\cC)\,.$$
Therefore
$$ \Add((\cA \amalg \cB)^{\sharp}_{pe} ,\cC) = \Add(\cA, \cC) \oplus \Add(\cB, \cC)\,,$$
which shows that $\Add$ has direct sums. 
\end{proof}
%\begin{remark}\label{rem:monAdd}
%Since the derived bi-functor of remark~\ref{rem:bi-functor} is bi-triangulated, the symmetric monoidal structure on $\Ho(\Spcat)^{tr}$ (\ref{prop:monoidal}) induces naturally a monoidal structure on $\Add$, making the natural functor
%$$ \Ho(\Spcat)^{tr} \too \Add$$
%symmetric monoidal.
%\end{remark}
%\begin{lemma}\label{lem1:T(A)s}
%Let $F:\Spcat \too \A$ be a functor with values in an additive category $\A$. Suppose that $F$ satisfies condition (M) of definition~\ref{def:additive}. Then $F$ satisfies condition (A), if and only if the induced functor $F: \Ho(\Spcat)^{tr} \to \A$ (\ref{thm:localization}) satisfies the following condition:
%\begin{itemize}
%\item[(A')] For every $\cB \in \Ho(\Spcat)^{tr}$, the inclusion morphisms of remark~\ref{rem:Tpe}, induce an isomorphism
%$$[F(i_1) \,\, F(i_2)] : F(\cB) \oplus F(\cB)\too F(T_{pe}(\cB))$$
%in $\A$.
%\end{itemize}
%\end{lemma}
%\begin{proof}
%By lemma~\ref{lem:T(A)s}, condition (A') implies condition (A). Now, let $\cB \in \Ho(\Spcat)^{tr}$. Since $\theta: \cB \stackrel{\sim}{\to} \cB^{\sharp}_{pe}$ is an isomorphism, lemma~\ref{lem:T(A)s} show us also that condition (A) implies condition (A').
%\end{proof}
\begin{notation}\label{def:rept}
Let $\cA$ and $\cB$ be two spectral categories, with $\cB$ fibrant. We denote by $\overline{\rep}(\cA, \cB)$ the full subcategory of $\cD(\cA \Dsma \cBo)$, whose objects are the $\cA \wedge \cBo$-modules $M$ such that $i_x^{\ast}(M) \in \cD_c(\cBo)$, for all objects $x \in \cA$ (see \ref{def:rep}). 
\end{notation}

\begin{remark}\label{rem:rept}
The category $\overline{\rep}(\cA,\cB)$ is naturally triangulated and we have a fully faithful functor $\rep(\cA,\cB) \to \overline{\rep}(\cA,\cB)$. Moreover, the morphism $\cB \to \Bpe$ induces an equivalence $\overline{\rep}(\cA,\cB) \stackrel{\sim}{\to} \rep(\cA,\Bpe)$ of triangulated categories.
\end{remark}
%\begin{proof}
%By lemma~\ref{lem:monoidal}, the morphism
%$$ \cA\op \Dsma \cB \simeq \cB \Dsma \cAo \too \Bpe \Dsma \cAo \simeq \cA\op \Dsma \Bpe$$
%is a Morita morphism. We obtain then an equivalence
%$$ \cD(\cA \Dsma \cBo) \stackrel{\sim}{\to} \cD(\cA \Dsma (\Bpe)\op)$$
%which restricts to the equivalence
%$$ \orep(\cA,\cB) \stackrel{\sim}{\too} \rep(\cA, \Bpe)\,.$$
%\end{proof}
\begin{proposition}\label{quasi-main}
Let $F: \Ho(\Spcat)^{tr} \to \A$ be a functor with values in an additive category $\A$. Then the following conditions are equivalent:
\begin{itemize}
\item[(1)] The functor $F$ is the composition of an additive functor $\Add \to \A$ with the natural functor $\Ho(\Spcat)^{tr} \to \Add$.
\item[(2)] For all $\cA, \cB \in \Ho(\Spcat)^{tr}$, the identity $F([X])+F([Z]) = F([Y])$ holds in $\A(F(\cA), F(\cB))$, for every triangle $X \to Y \to Z \to \Sigma X$ in $\rep(\cA,\cB)$.
\item[(3)] For every $\cB \in \Ho(\Spcat)^{tr}$, the induced map 
$$[F((i_1)^{\sharp}_{pe}) \,\, F((i_2)^{\sharp}_{pe})] : F(\Bpe) \oplus F(\Bpe)\too F(T(\cB)^{\sharp}_{pe})$$
is an isomorphism in $\A$.
\end{itemize} 
\end{proposition}
\begin{proof}
By construction of $\Add$, conditions (1) and (2) are equivalent. We now show that condition (1) implies condition (3). Since by hypothesis $F$ factors through $\Add$, it is enough to show that the induced map
$$[(i_1)^{\sharp}_{pe} \,\, (i_2)^{\sharp}_{pe}] : \Bpe \oplus \Bpe \too T(\cB)^{\sharp}_{pe}$$
is an isomorphism in $\Add$. By the Yoneda Lemma, we need to show that for every $\cC \in \Add$, the induced map (on the top)
$$ 
\xymatrix{
\Add(\cC, \Bpe \oplus \Bpe)  \ar[r] & \Add(\cC, T(\cB)^{\sharp}_{pe}) \ar@{=}[d] \\
K_0\rep(\cC, \Bpe) \oplus K_0\rep(\cC, \Bpe) \ar[u]^{\sim} \ar[r] & K_0\rep(\cC, T(\cB)^{\sharp}_{pe})
}
$$
is an isomorphism. By Remark~\ref{rem:rept}, we have the following equivalences
$$
\begin{array}{lcr}
\overline{\rep}(\cC, \cB) \simeq \rep(\cC, \Bpe) && \overline{\rep}(\cC, T(\cB)) \simeq \rep(\cC, T(\cB)^{\sharp}_{pe})\,.
\end{array}
$$
Therefore, it is enough to show that the triangulated category $\orep(\cC, T(\cB))$ admits a semi-orthogonal decomposition~\cite[2.4]{Bon-Orl} in two subcategories equivalent to $\orep(\cC,\cB)$. Given an object $X \in \orep(\cC, T(\cB))$, we can represent it by a cofibration $X_1 \to X_2$ between fibrant and cofibrant objects in $\orep(\cC,\cB)$. The triangle associated to $X$ by the semi-orthogonal decomposition is then the one induced by the following diagram
$$
\xymatrix{
X_1 \ar@{=}[r] & *+<3ex>{X_1} \ar[r] \ar@{>->}[d]& 0 \ar[d]\\
*+<3ex>{X_1}  \ar@{=}[u] \ar@{>->}[r] & X_2  \ar[r] & X_2/X_1 \,,
}
$$
where the left and right vertical terms belong to $\orep(\cC,\cB)$.
We now show that condition (3) implies condition (2). Let $$X \to Y \to Z \to \Sigma X$$ be a triangle in $\rep(\cA,\cB)$. It is isomorphic in $\rep(\cA,\cB)$ to a triangle associated with a homotopy cofiber sequence
$$ X' \too Y' \too Z'$$
of $\cA\wedge\cBo$-modules. We can then consider $M:=(X' \to Y')$ as an element of $\rep(\cA, T(\cB))$. By Theorem~\ref{thm:Hom}, the isomorphism class $[M]$ of $M$ in $\rep(\cA, T(\cB))$ corresponds to a morphism in $\Ho(\Spcat)$ from $\cA$ to $T(\cB)$. By composing it with $\theta: T(\cB) \to T(\cB)^{\sharp}_{pe}$, we obtain a morphism $\theta\circ M$ in $\Ho(\Spcat)^{tr}$ from $\cA$ to $T(\cB)^{\sharp}_{pe}$. Now, since $\cB \in \Ho(\Spcat)^{tr}$, we have by Remark~\ref{rk:adjoints1} induced morphisms 
$$
\begin{array}{lcccr}
\xymatrix{
T(\cB)^{\sharp}_{pe} \ar@<-6ex>[d]_{(P_1)^{\sharp}_{pe}} \ar[d]^{(P)^{\sharp}_{pe}} \ar@<7ex>[d]^{(P_2)^{\sharp}_{pe}} \\
\Bpe
} & & & & 
\xymatrix{
T(\cB)^{\sharp}_{pe} \\
\Bpe \ar@<-1ex>[u]_{(i_2)^{\sharp}_{pe}} \ar@<1ex>[u]^{(i_1)^{\sharp}_{pe}}
}
\end{array}
$$
in $\Ho(\Spcat)^{tr}$, satisfying the equalities
$$
\begin{array}{cccc}
(P_1 \circ i_1)^{\sharp}_{pe} = \mbox{Id} & (P_1 \circ i_2)^{\sharp}_{pe} =\ast & (P_2\circ i_1)^{\sharp}_{pe} =\ast & (P_2 \circ i_2)^{\sharp}_{pe}=\mbox{Id}
\end{array}
$$
and
$$
\begin{array}{cc}
(P \circ i_1)^{\sharp}_{pe}=\mbox{Id} & (P \circ i_2)^{\sharp}_{pe}=\mbox{Id}\,.
\end{array}
$$
This implies that in the additive category $\A$, we have the equalities
$$(F((P_1)^{\sharp}_{pe}) + F((P_2)^{\sharp}_{pe}) \circ [F((i_1)^{\sharp}_{pe})\,\,\,F((i_2)^{\sharp}_{pe})] = [\mbox{Id} \,\,\, \mbox{Id}]$$
and
$$ F((P)^{\sharp}_{pe}) \circ [F((i_1)^{\sharp}_{pe}) \,\,\, F((i_2)^{\sharp}_{pe})] =[\mbox{Id} \,\,\, \mbox{Id}]\,.$$
By hypothesis, the morphism $[F((i_1)^{\sharp}_{pe})\,\,\, F((i_2)^{\sharp}_{pe})]$ is invertible and so we obtain the following equality
$$ F((P_1)^{\sharp}_{pe}) + F((P_2)^{\sharp}_{pe}) = F((P)^{\sharp}_{pe})\,.$$
Since $\cB \in \Ho(\Spcat)^{tr}$, the morphism $\theta: \cB \stackrel{\sim}{\to} \Bpe$ is an isomorphism and so $\rep(\cA, \cB) \simeq \rep(\cA,\Bpe)$. Using the extensions of scalars functors of Remark~\ref{rk:adjoints1}, we observe that
$$ 
\begin{array}{ccc}
(P_2)^{\sharp}_{pe}\circ [\theta \circ M]=[X] & (P)^{\sharp}_{pe} \circ [\theta \circ M]=[Y] & (P_1)^{\sharp}_{pe} \circ [\theta \circ M]=[Z]\,,
\end{array}
$$
and so we conclude that the identity
$$ F([X]) +F([Z]) = F([Y])$$
holds in $\A(F(\cA),F(\cB))$.
\end{proof}
\begin{remark}\label{rk:EasyAdd}
Notice that the argument used in the proof of the implication $(1) \Rightarrow (3)$ in Theorem~\ref{quasi-main}, shows us that for any upper triangular matrix $\underline{M}$, the induced map
$$[(i_1)^{\sharp}_{pe} \,\, (i_2)^{\sharp}_{pe}] : \cA^{\sharp}_{pe} \oplus \cC^{\sharp}_{pe} \too |\underline{M}|^{\sharp}_{pe}$$
is an isomorphism in $\Add$. This implies that a functor $F: \Spcat \to \A$, which satisfies condition M), satisfies condition MA) if and only if it satisfies the following condition
\begin{itemize}
\item[MA')] For every $\cB \in \Ho(\Spcat)^{tr}$, the induced map 
$$ [F(i_1) \,\, F(i_2)]: F(\cB) \oplus F(\cB) \too F(T(\cB))$$
is an isomorphism in $\A$.
\end{itemize}
\end{remark}

\subsection{Main Theorem}\label{ssec:univ}
\begin{theorem}\label{thm:main}
The composed functor
$$\cU: \Spcat \too \Ho(\Spcat) \stackrel{(-)^{\sharp}_{pe}}{\too} \Ho(\Spcat)^{tr} \too \Add$$
is the universal matrix invariant of spectral categories, \ie for every additive category $\A$, the functor $\cU$ induces a bijection between the additive functors from $\Add$ to $\A$ and the matrix invariants of spectral categories (\ref{def:additive}). 
\end{theorem}
\begin{proof}
By Proposition~\ref{quasi-main}, the composition of $\cU$ with an additive functor from $\Add$ to $\A$ gives rise to an matrix invariant of spectral categories. Now, let $F:\Spcat \to \A$ be a matrix invariant. We must show that $F$ factors uniquely through $\cU$. Since $F$ satisfies condition M), Theorem~\ref{thm:localization} implies that $F$ factors uniquely through the composition
$$ \Spcat \too \Ho(\Spcat) \too \Ho(\Spcat)^{tr}\,.$$
Since $F$ satisfies condition MA), the induced functor
$ \Ho(\Spcat)^{tr} \to \A$
satisfies moreover the condition (3) of Theorem~\ref{quasi-main}, and so the proof is finished.
\end{proof}

%-----------------------------------------------------------------------
\section{Algebraic $K$-theory}\label{sec:AlgK}
%-----------------------------------------------------------------------
Let $\cA$ be a spectral category. We denote by $\per(\cA) \subset \Aomod$ the category of cofibrant and {\em perfect} (\ie compact in $\cD(\cAo)$) $\cAo$-modules. The category $\per(\cA)$ carries a natural Waldhausen structure, in the sense of~\cite[\S\,3]{DS}. The {\em $K$-theory spectrum $K(\cA)$ of $\cA$} is defined by applying the Waldhausen's $S_{\bullet}$-construction \cite{Wald} to $\per(\cA)$. Let $F:\cA \to \cB$ be a spectral functor. Since the extension of scalars functor $F_!:\Aomod \to \Bomod$ preserves pushouts and cofibrant and perfect objects, it induces a Waldhausen's functor $F_!:\per(\cA) \to \per(\cB)$. In sum, we obtain a well-defined functor
$$ K(-): \Spcat \too \Ho(\Spe)\,,$$
with values in the homotopy category of spectra.
\begin{proposition}\label{prop:AddKth}
The algebraic $K$-theory functor
$$ K(-): \Spcat \too \Ho(\Spe)$$
is a matrix invariant of spectral categories (\ref{def:additive}).
\end{proposition}
\begin{proof}
Let $F:\cA \to \cB$ be a Morita equivalence. Since $ \bbL F_!: \cD(\cAo) \stackrel{\sim}{\to} \cD(\cBo)$
is an equivalence, \cite[Prop.\,3.7]{DS} implies that $$K(F):K(\cA) \stackrel{\sim}{\too} K(\cB)$$ is an isomorphism in $\Ho(\Spe)$. Condition MA) follows from~\cite[Thm.\,1.4]{Wald}.
\end{proof}
%----------------------------------------------------------------------------------------
\section{Topological Hochschild homology}\label{sec:THH}
%----------------------------------------------------------------------------------------
Recall from \cite[\S\,3]{Mandell} \cite[\S 10]{SpectralAlg} the construction of the topological Hochschild homology ($THH$) and topological cyclic homology ($TC$) of spectral categories. We have well-defined functors
$$ THH(-), \,TC(-): \Spcat \too \Ho(\Spe)\,.$$
\begin{proposition}\label{prop:THHAdd}
The functors $THH(-)$ and $TC(-)$ are matrix invariants of spectral categories.
\end{proposition}
\begin{proof}
By \cite[Prop.\,3.8-3.9]{Mandell}, it is enough to show that $THH(-)$ is a matrix invariant. Let us start by showing that $THH(-)$ satisfies condition M). By \cite[Thm.\,4.9]{Mandell}, $THH(-)$ sends the stable quasi-equivalences to isomorphisms. Therefore, $THH(-)$ descends to $\Ho(\Spcat)$. An application of \cite[Thm.\,4.12]{Mandell} (with $\cC=\cA$, $\cC'=\Ape$ and $\cD= \Aomod$), shows us that for every $\cA \in \Ho(\Spcat)$, the morphism
$$ \theta:\cA \too \Ape\,,$$
is sent to an isomorphism by $THH(-)$. This implies that $THH(-)$ sends the Morita equivalences to isomorphisms. We now show that $THH(-)$ satisfies condition MA') of remark~\ref{rk:EasyAdd}. Let $\cB$ be an object in $\Ho(\Spcat)^{tr}$. Consider the following diagram of spectral categories (see~\ref{rk:adjoints1})
$$ 
\xymatrix{
\cB \ar[r]_{i_1} & T(\cB)  \ar@<-1ex>[l]_{P_1}
\ar[r]_{P_2} & \cB \ar@<-1ex>[l]_{i_2} \,.
}
$$
This diagram gives rise to an exact sequence of triangulated categories
$$ 0 \too \cD(\cBo) \too \cD(T(\cBo)) \too \cD(\cBo) \too 0\,,$$
and so by \cite[Thm.\,6.1]{Mandell} to a distinguished triangle
$$ THH(\cB) \too THH(T(\cB)) \too THH(\cB) \too THH(\cB)[1]$$
in $\Ho(\Spe)$. Since $P_2 \circ i_2 = \mbox{Id}$ and $P_1 \circ i_1 = \mbox{Id}$, this triangle splits and so the induced map
$$[THH(i_1) \,\,THH(i_2)]: THH(\cB) \oplus THH(\cB) \too THH(T(\cB))$$
is an isomorphism in $\Ho(\Spe)$.
\end{proof}
%-------------------------------------------------------------------------------
\section{Grothendieck group and trace maps}\label{sec:Grothendieck}
%-------------------------------------------------------------------------------
Let $\underline{\bbS}$ be the spectral category with one object and endomorphisms symmetric ring spectrum the sphere symmetric ring spectrum $\bbS$. We denote by $K_0(\cB)$ be the Grothendieck group of a spectral category $\cB$, \ie the $0$-th stable homotopy group of $K(\cB)$ (see Chapter \ref{sec:AlgK}) or equivalently the Grothendieck group of the triangulated category $\cD_c(\cBo)$. 
\begin{proposition}\label{prop:co-repres}
For every spectral category $\cB$, we have a natural isomorphism of abelian groups
$$ \Add(\cU(\underline{\bbS}), \cU(\cB)) \stackrel{\sim}{\too} K_0(\cB)\,.$$
\end{proposition}
\begin{proof}
Notice that we have the following isomorphisms
\begin{eqnarray}
 \Add(\cU(\underline{\bbS}), \cU(\cB)) & = & K_0\rep(\underline{\bbS}^{\sharp}_{pe}, \cB^{\sharp}_{pe}) \nonumber \\
 & \simeq & K_0\rep(\bbS, \Bpe) \\
 & \simeq & K_0\orep(\bbS,\cB) \\
 & \simeq & K_0\cD_c(\cBo)\,, \nonumber
\end{eqnarray}
where (11.0.5) is follows from Remark~\ref{rk:eqrep} and (11.0.6) from Remark~\ref{rem:rept}. 
\end{proof}
This co-representability result has the following important application: 
\begin{corollary}\label{cor:trace}
Let $j$ be a non-negative integer and
$$ 
\begin{array}{lcr}
K_0(-):\Spcat \too \Ab &  & THH_j(-) : \Spcat \too \Ab\,,
\end{array}
$$
the Grothendieck and the $j$-th topological Hochschild homology group functors. Then, each generator $g$ of the $j$-th stable homotopy group of the sphere, furnishes for free a non-trivial trace map
$$ tr_{j,g}:K_0(-) \Rightarrow THH_j(-)\,.$$
\end{corollary}
\begin{proof}
Since $K_0(-)$ and $THH_j(-)$ are matrix invariants, they descend by Theorem~\ref{thm:main}, to two additive functors $\overline{K_0}(-)$ and $\overline{THH_j}(-)$ defined on $\Add$. By Proposition~\ref{prop:co-repres} and the (enriched) Yoneda Lemma, we have a natural isomorphism
$$ \mbox{Nat}(\overline{K_0}(-), \overline{THH_j}(-)) \simeq \overline{THH_j}(\cU(\underline{\bbS})) = THH_j(\bbS)\,,$$
where $\mbox{Nat}(-,-)$ denotes the abelian group of natural transformations. Since $THH_j(\bbS) = \pi^s_j(\bbS)$, each generator $g$ of $\pi_j^s(\bbS)$ furnishes us a natural transformation
$$ \overline{tr_{j,g}}:\overline{K_0}(-) \Rightarrow \overline{THH_j}(-)\,,$$
whose pre-composition with $\cU$ gives rise to a non-trivial trace map
$$ tr_{j,g}:K_0(-) \Rightarrow THH_j(-)\,.$$
\end{proof}


\begin{thebibliography}{00}

%\bibitem{Bergner} J.~Bergner, {\em A Quillen model structure on the category of simplicial categories}Trans.~Amer.~Math.~Soc. {\bf 359} (2007), 2043--2058.

\bibitem{Grothendieck}  M.~Artin, A.~Grothendieck, J. L.~Verdier, {\em Theorie des topos et cohomologie {\'e}tale des
sch{\'e}mas}. SGA4, Tome 3, Lecture Notes in Mathematics, vol. {\bf 305}, Springer Verlag, 1973.  

\bibitem{Balmer} P.~Balmer, M.~Schlichting, {\em Idempotent completion
    of triangulated categories}, J.~Alg. {\bf 236} (2001), no.~2, 819--834. 

\bibitem{Mandell} A.~Blumberg, M.~Mandell, {\em Localization theorems in topological Hochschild homology and topological cyclic homology}. Available at arXiv:$0802.3938$.    

\bibitem{Bon-Orl} A.~Bondal, D.~Orlov, {\em Semiorthogonal
    decomposition for algebraic varieties}, preprint MPIM 95/15
    (1995), preprint math. AG/9506012.

\bibitem{Cortinas} G.~Corti{\~n}as, A.~Thom, {\em Bivariant algebraic $K$-theory}, J.~Reine Angew. Math, {\bf 510}, 71--124.

%\bibitem{BK} A.~Bondal, M.~Kapranov, {\em Framed triangulated categories} (Russian) Mat.~Sb. {\bf 181} (1990) no.~5, 669--683; translation in Math.~USSR-Sb. {\bf 70} no.~1, 93--107.

%\bibitem{Borceaux} F.~Borceux, {\em Handbook of categorical algebra. 2}. Encyclopedia of Mathematics and its Applications, {\bf 51}, Cambridge University Press, 1994. 

%\bibitem{Drinfeld} V.~Drinfeld, {\em DG quotients of DG categories}, J. Algebra {\bf 272} (2004), 643--691.

%\bibitem{Chitalk} V.~Drinfeld, {\em DG categories}. University of Chicago Geometric Langlands Seminar. Notes available at {\tt http://www.math.utexas.edu/users/benzvi/GRASP/lectures/Langlands.html}.

\bibitem{Dugger} D.~Dugger, {\em Spectral enrichments of model categories}, Homology, Homotopy, Appl. {\bf 8} (2006), no. 1, 1-30.

\bibitem{DS} D.~Dugger, B.~Shipley, {\em $K$-theory and derived equivalences}, Duke Math. J. {\bf 124} (2004), no.3, 587--617.

 \bibitem{EKMM} A.~Elmendorf, I.~Kriz, M.~Mandell, P.~May, {\em Rings, Modules, and Algebras in Stable Homotopy Theory}, with an appendix by M. Cole. Mathematical Surveys and Monographs, {\bf47}. American Mathematical Society, Providence, RI, 1997.

%\bibitem{Dwyer}, {\em Homotopy theories}. Handbook of algebraic topology.....

\bibitem{Zisman} P.~Gabriel, M.~Zisman, {\em Calculus of fractions and homotopy theory}, Ergebnisse der Mathematik und ihrer Grenzgebiete, Band {\bf 35}, Springer-Verlag New York, Inc., New York, 1967.

\bibitem{Garkusha} G.~Garkusha, {\em Homotopy theory of associative rings}, Adv. in Math., {\bf 213}, no.2, 553-599.

\bibitem{GH} T.~Geisser, L.~Hesselholt, {\em On the $K$-theory and topological cyclic homology of smooth schemes over a discrete valuation ring}. Trans.~Amer.~Math.~Soc. {\bf 358}(1), 131--145, 2006.

%\bibitem{GJ} P. Goerss and J. Jardine, {\em Simplicial homotopy theory}, Progress in Mathematics 174, Birkh{\"a}user.

\bibitem{HSS} M.~Hovey, B.~Shipley, J.~Smith, {\em Symmetric spectra}. J. Amer. Math. Soc. {\bf 13} (2000), no. 1, 149--208.

\bibitem{Hirschhorn} P.~Hirschhorn, {\em Model categories and their localizations}, Mathematical Surveys and Monographs, {\bf 99}, American Mathematical Society, 2003.
    
% \bibitem{Hovey} M.~Hovey, {\em Model categories}, Mathematical Surveys and Monographs, {\bf 63}, American Mathematical Society, 1999.
    
%\bibitem{ENS} M.~Kontsevich, {\em  Categorification, NC Motives, Geometric Langlands and Lattice Models}. University of Chicago Geometric Langlands Seminar. Notes available at {\tt http://www.math.utexas.edu/users/benzvi/notes.html}. 

\bibitem{ICM} B.~Keller, {\em On differential graded categories}. International Congress of Mathematicians, Vol.~II,  151--190, Eur.~Math.~Soc., Z{\"u}rich, 2006.

\bibitem{Kontsevich} M.~Kontsevich, {\em Non-commutative motives}. Talk at IAS, October 2005. Video available at {\tt http://video.ias.edu/Geometry-and-Arithmetic}.

\bibitem{finMotiv} \bysame, {\em Notes on motives in finite characteristic}. Preprint arXiv:$0702206$. To appear in Manin Festschrift.

\bibitem{Lurie} J.~Lurie, {\em Higher topos theory}. Available at arXiv:$0608040$. 

\bibitem{Lydakis} M.~Lydakis , {\em Simplicial functors and stable homotopy theory}. Preprint (1998). Available at {\text http://hopf.math.purdue.edu}.

\bibitem{May} M.~A.~Mandell, J.~P.~May, S.~Schwede, B.~Shipley, {\em  Model categories of diagram spectra}. Proc. London Math. Soc. (3) {\bf 82} (2001), no. 2, 441--512.    
        
\bibitem{MacLane}  S.~Mac Lane, {\em Categories for the working
    mathematician}, Graduate Texts in Mathematics, {\bf 5},
  Springer-Verlag, 1998.

\bibitem{Nest} R.~Meyer, R.~Nest, {\em The Baum-Connes conjecture via localisation of categories}, Topology, {\bf 45}, 2006, 209--259.

\bibitem{Neeman} A.~Neeman, {\em Triangulated categories}, Annals of Mathematics Studies, {\bf 148}, Princeton University Press, 2001.

\bibitem{Quillen} D.~Quillen, {\em Homotopical algebra}, Lecture Notes
  in Mathematics, {\bf 43}, Springer-Verlag, 1967.

\bibitem{Schwede} S.~Schwede, {\em An untitled book project about symmetric spectra}. Available at {\tt www.math.uni-bonn.de/people/schwede}.

%\bibitem{SS1} S.~Schwede, B.~Shipley, {\em Algebras and modules in monoidal model categories}. Proc. London Math. Soc. (3) {\bf 80} (2000), no. 2, 491--511.

\bibitem{SS} S.~Schwede, B.~Shipley, {\em Stable model categories are categories of modules}. Topology {\bf 42} (2003), no. 1, 103--153.

\bibitem{Spectral} G.~Tabuada, {\em Homotopy theory of Spectral categories}. Available at arXiv:$0801.4524$(version~3). To appear in Adv. in Math.

\bibitem{IMRN} \bysame, {\em Invariants additifs de dg-cat{\'e}gories}, Int.~Math.~Res.~Notices {\bf 53} (2005), 3309--3339.

\bibitem{IMRNC} \bysame, {\em Corrections {\`a} Invariants additifs de dg-cat{\'e}goriesÕ}, Int.~Math.~Res.~Notices, article ID rnm {\bf 149} (2007).

\bibitem{Duke} \bysame, {\em Higher $K$-theory via universal invariants}. Duke Math. J. {\bf 145} (2008), no.1, 121--206.

\bibitem{SpectralAlg} \bysame, {\em Topological Hochschild and cyclic homology for Differential graded categories}. Available at arXiv:$0804.2791$.

%\bibitem{TV} B.~To{\"e}n and M.~Vaqui{\'e}, {\em Moduli of objects in dg-categories}. Available at arXiv:math/0503269. To appear in Annales de l'E.N.S.

\bibitem{Toen} B.~To{\"e}n, {\em The homotopy theory of dg-categories and
    derived Morita theory}, Invent. Math. 167 (2007), no.~3, 615--667.
    
\bibitem{Toen1} \bysame, {\em Lectures on dg-categories}. Sedano winter school on K-theory, January 2007. Available at {\tt http://www.math.univ-toulouse.fr/toen/note.html}.    

\bibitem{Voevodsky} V.~Voevodsky, {\em Triangulated categories of motives over a field}. Cycles, transfers, and
motivic homology theories, Annals of Mathematics Studies, vol. {\bf 143}, Princeton University
Press, Princeton, NJ, 2000, 188Ð238.

\bibitem{Wald} F.~Waldhausen, {\em Algebraic K-theory of spaces},
  Algebraic and geometric topology (New Brunswick, N.~J., 1983),
  318--419, Lecture Notes in Math., 1126, Springer, Berlin, 1985.

\end{thebibliography}
\end{document}